# UTILIZING THE REDUNDANT CONSTRAINTS FOR THE UPLIFT PAYMENT ELIMINATION

Vadim Borokhov[1]

LLC "En+development",
ul. Vasilisi Kozhinoy 1, Moscow, 121096, Russia

**Abstract**
A power market with non-convexities may not have an equilibrium price for power that provides economic stability of the centralized dispatch outcome. In this case, the market players are entitled to receive the uplift payments that compensate the economic profit lost when following the centralized dispatch. We consider a special class of the (possibly non-linear) redundant constraints that are redundant not only on the feasible set of the centralized dispatch optimization problem (and, therefore, do not change the centralized dispatch outcome) but also on the larger set obtained when the power balance constraint is relaxed. We show that the Lagrangian relaxation of these redundant constraints may reduce the uplift payments without changing the duality gap. For any given market price (or a pricing algorithm that sets the producer revenue as a function of its output volume) in a uninode multi-period power market with fixed load, we explicitly construct a family of the redundant constraints that do not change the maximum profit of the producer and result in zero uplift payment. We show that the introduction and subsequent Lagrangian relaxation of just one redundant constraint in the centralized dispatch problem suffice to eliminate the uplift payments for all the producers. In the case of the convex hull pricing method, the introduction of these redundant constraints affects neither the duality gap nor the market price for power. The results can be straightforwardly generalized to a power market with the price-sensitive load.

## I. Introduction

Many deregulated electric power markets are centrally coordinated with generating unit dispatch and load schedule obtained from a bid-based security-constrained centralized dispatch optimization problem. The solution of this problem also produces the system marginal price (or locational marginal prices) for power [1]-[3]. If the optimization problem is convex, then the marginal price is an equilibrium price for power and no market player (acting as a price-taker) has the economic incentives to deviate from the centralized dispatch outcome. However, if the optimization problem is not convex, then the marginal price may fail to ensure the economic equilibrium of the centralized dispatch outcome since the non-convex components of the generator and consumer bids are not affecting the value of the marginal price. The non-convexities usually originate both from the supply side (due to the no-load cost, start-up cost, non-zero minimum output limits, integral commitment decision variables, minimum up/down times, etc.) and from the flexible demand side due to discrete and minimum power consumption levels [5]. In fact, in case of the non-convex centralized dispatch optimization problem, an equilibrium market price may not exist at all and some other pricing approach has to be implemented to ensure the economic stability of the centralized dispatch outcome [4]. Recently, a number of pricing methods for the centrally coordinated markets with non-convexities have been developed, [6]-[24]. The new service (a unit being online) and the corresponding unit specific prices were introduced in [6] by constraining the integral status variables to their optimal values obtained from the centralized dispatch. However, these new prices can be negative for some units and the method is similar to pay-as-bid pricing. If only the positive prices are applied to allow generators to retain their profits, then the competitive equilibrium is not achieved. In [7], [8] this method was further improved to generate more stable prices by adding extra constraints to the reformulated optimization problem that also set certain continuous variables to their optimal values. The nonlinear (discriminatory) pricing methods for power with market player specific prices were developed in [9]-[12]. The nonlinear pricing in the form of the generalized uplift functions that includes generators as well as consumers in the lost profit compensation and ensures zero net uplift at the market was proposed in [13]-[15]. The minimum zero-sum uplift pricing approach that increases the price above marginal cost and transfers all the additional payments (that the profitable suppliers receive as a result of the price increase) to the unprofitable suppliers to make them whole in the form of internal zero-sum uplifts was introduced in [16]. In [17] a primal-dual approach was proposed to find the market prices that minimize the social welfare reduction due to schedules inconsistency and ensure non-negative generator profits. However, in this approach, some of the lost profit may not be compensated to generators and the competitive equilibrium at the centralized dispatch solution is not achieved. In [18] a semi-Lagrangian relaxation approach was developed to find a uniform market price that produces the same solution as the original centralized dispatch problem while ensuring the non-confiscatory pricing for generators. A zero-sum uplift pricing scheme that minimizes the maximum contribution to the uplift financing in a market with price-sensitive load was suggested in [19]. In the case of no price-sensitive load, this approach produces the market price equal the maximum average cost of the generators. The minimum-uplift pricing (also known as the convex hull pricing) was proposed in [20]-[22] and yields a uniform market price that minimizes the total uplift payment needed to ensure the economic stability of the centralized dispatch outcome. In this approach, at a given market price each market player is compensated the lost profit calculated as the difference between the maximum value of its profit function on the market player private feasible set and its profit received when following the centralized dispatch. Since the uplift payments distort the uniform market pricing and decrease the transparency of the market pricing method, it is critical to reduce these payments. In [23], [24] it was proposed to modify the minimum-uplift pricing method by excluding the power volumes that are not attainable in a decentralized market from the lost profit calculation since the opportunities to supply these volumes are not forgone by a market player when accepting the centralized dispatch outcome. This approach results in the lower (or equal) total uplift payment compared to the minimum-uplift pricing algorithm.

---

[1] E-mail: vadimbor@yahoo.com. The views expressed in this paper are solely those of the author and not necessarily those of LLC "En+development".

For the convex hull pricing method, in [25] it was suggested to reduce the total uplift payment, which compensates the lost profit of the market players, at the expense of having one affine redundant constraint introduced in the centralized dispatch optimization problem. This new constraint depends on the unit status variables of all generators and leads to the introduction of the new service (a unit being online) and the associated price in addition to the market price for power, which can be viewed as the producer revenue function amendments. The linearity of the redundant constraint ensures that the duality gap, introduced by the Lagrangian relaxation of both the power balance constraint and the new constraint, is equal to that in the absence of the redundant constraint [26]. However, introduction and subsequent dualization of the new constraint entails that the duality gap may no longer coincide with the total uplift, which is potentially reduced but generally still non-zero.

In this paper, we study the problem of the total uplift (lost profit) reduction in a general pricing setting, which fixes the producer (consumer) revenue (cost) as a function of its status-output (consumption) variables, by introduction of the redundant constraints and the corresponding non-negative amendments to the revenue (cost) functions. Thus, our study is also applicable to the cases with uniform pricing for power (such as marginal pricing, convex hull pricing) and discriminatory (non-linear) pricing with the uplift payments. For simplicity, we consider a multi-period uninode power market with the fixed load. The analysis and the results can be easily translated to markets with price-sensitive demand. To simplify the notations, we assume that each producer operates just one generating unit.

We consider a special type of the redundant constraints – the constraints that hold not only on the feasible set of the centralized dispatch optimization problem but also on the larger set obtained by relaxing the power balance constraint. The redundant constraints under consideration are introduced in the market player individual profit optimization problems. Therefore, we require that each constraint depends on just one producer status-output variable. We show that it suffices to consider only this class of the redundant constraints to fully absorb the uplift payment of a producer (thus, resulting in zero uplift payment) and find the general form expression for the corresponding revenue amendment function for the producer. This function satisfies the following three properties: it is non-negative on the producer private feasible set, makes no contribution to the producer maximum profit, and yields zero uplift payment for the producer. For the uniform market price, we show that just one redundant constraint, which is the sum (over all the producers) of these properly rescaled redundant constraints) introduced directly in the centralized dispatch optimization problem suffices to produce zero total uplift. If the uniform market price for power is set by the convex hull pricing method, the dualization of both the power balance constraint and the new (redundant) constraint results in the same set of the market prices and the same maximum profit for each producer but gives zero total uplift.

The paper is organized as follows. In Section 2 we formulate the conditions on the revenue amendment function that is non-negative on the producer private feasible set, leaves the producer maximum profit unaffected, and fully absorbs the uplift payment. In Section 3 we introduce the redundant constraints and study their relations with the revenue amendment functions and the uplift reduction problem. In Section 4 we formulate the necessary and sufficient conditions for a given set of the redundant constraints and the associated multipliers to produce zero uplift. The general form expression for the revenue amendment function that satisfies the three abovementioned properties is obtained in Section 5. In Section 6 we apply the proposal to power markets with marginal pricing, while the application of the method to a producer with linear cost function is given in Section 7. In Section 8 we construct the revenue amendment functions and the corresponding redundant constraints in a numerical example. The results are summarized in Section 9. Some mathematical aspects of the redundant constraints are summarized in the Appendix.

## II. The problem formulation

Consider $T$-period uninode power market with the fixed load $\mathbf{d}$, $\mathbf{d} \in R_{\geq 0}^T$, with $|I|$ producers, where $|\bullet|$ denotes cardinality of a set. For each producer $i \in I$ at time period $t \in \{1,...,T\}$, let $u_i^t$ and $g_i^t$ denote the commitment and output variables, respectively. Introduce $\mathbf{g}_i = (g_i^1,...,g_i^T)$, $\mathbf{u}_i = (u_i^1,...,u_i^T)$, $x_i^t = (u_i^t, g_i^t)$, $\mathbf{x}_i = (x_i^1,...,x_i^T)$, $\mathbf{x} = (\mathbf{x}_1, \ldots, \mathbf{x}_{|I|})$. Let $X_i$ be the producer $i$ private feasible set (which is assumed to be nonempty and compact) and $C_i(\mathbf{x}_i)$ be the offer cost function of the producer $i$. The centralized dispatch problem has the form

$$f^* = \min_{\substack{\mathbf{x}_i \in X_i, \forall i \in I \\ \sum_{i \in I} \mathbf{g}_i = \mathbf{d}}} \sum_{i \in I} C_i(\mathbf{x}_i) \qquad (1)$$

The feasible set of (1) is assumed to be nonempty and compact. Let $\mathbf{x}^* = (\mathbf{x}_1^*, \ldots, \mathbf{x}_{|I|}^*)$ denote an optimal point of (1). Although we consider a centrally coordinated power market with a fixed load, it is straightforward to include the price-sensitive demand in our analyses. For a given price $\mathbf{p}$, $\mathbf{p} \in R^T$, the standard expression for the producer $i$ profit function is given by $\pi_i^{st.}(\mathbf{p}, \mathbf{x}_i) = R_i^{st.}(\mathbf{p}, \mathbf{x}_i) - C_i(\mathbf{x}_i)$ with the revenue function $R_i^{st.}(\mathbf{p}, \mathbf{x}_i)$ usually expressed as $R_i^{st.}(\mathbf{p}, \mathbf{x}_i) = \mathbf{p}^T \mathbf{g}_i$. However, our analysis is performed for the most general form[2] of the function $R_i^{st.}(\mathbf{p}, \mathbf{x}_i)$, unless explicitly stated otherwise. The generator status-output value determined by the centralized dispatch optimization problem results in the profit $\pi_i^{st.*}(\mathbf{p}) = R_i^{st.}(\mathbf{p}, \mathbf{x}_i^*) - C_i(\mathbf{x}_i^*)$. Given the price $\mathbf{p}$, the maximum producer's profit is $\pi_i^{st.+}(\mathbf{p}) = \max_{\mathbf{x}_i \in X_i} \pi_i^{st.}(\mathbf{p}, \mathbf{x}_i)$. If $\pi_i^{st.+}(\mathbf{p}) > \pi_i^{st.*}(\mathbf{p})$, then the producer has economic incentives to adjust its output volumes to attain the maximum profit and deviate from the centralized dispatch outcome. The standard procedure utilized to eliminate such incentives is to pay the producer the uplift in the amount of $\pi_i^{st.+}(\mathbf{p}) - \pi_i^{st.*}(\mathbf{p})$ if it follows the centralized dispatch outcome within a set tolerance band. In this case, the producer receives the profit $\pi_i^{st.+}(\mathbf{p})$ when supplying the output $\mathbf{g}_i^*$. (We make the usual assumption that a generator decides to deviate from $\mathbf{x}_i^*$ only if it receives a *higher* profit when having a different output volume. This means that if both $\mathbf{x}_i^*$ and $\mathbf{x}_i'$ maximize the generator profit, then it will not deviate from $\mathbf{x}_i^*$ to $\mathbf{x}_i'$.) The uplift payment $\pi_i^{st.+}(\mathbf{p}) - \pi_i^{st.*}(\mathbf{p})$ can be viewed as the cost of the commitment ticket payable to the generator $i$ for following the centralized dispatch [20]-[22]. If the price-sensitive demand is present in the system, then such a compensation mechanism should be applied to the demand side as well. (We note that for some $R_i^{st.}(\mathbf{p}, \mathbf{x}_i)$, in particular $R_i^{st.}(\mathbf{p}, \mathbf{x}_i) = \mathbf{p}^T \mathbf{g}_i$, if the total uplift payment is non-zero and all the consumers submit only the price-sensitive bids, then this leads to the budget-balancing problem as the total uplift payment (if non-zero) exceeds the amount that can be collected from the market players provided that no consumer (producer) can be charged (paid) above (below) its bid cost. In the present paper, we do not address this problem.) The uplift payment in the amount of $\pi_i^{st.+}(\mathbf{p}) - \pi_i^{st.*}(\mathbf{p})$ can be expressed as the generator revenue function amendment of the form $\delta_{\mathbf{x}_i, \mathbf{x}_i^*}[\pi_i^{st.+}(\mathbf{p}) - \pi_i^{st.*}(\mathbf{p})]$ with a function $\delta_{\mathbf{x}_i, \mathbf{x}_i^*}$ defined as $\delta_{\mathbf{x}_i, \mathbf{x}_i^*} = 1$ for $\mathbf{x}_i = \mathbf{x}_i^*$, and $\delta_{\mathbf{x}_i, \mathbf{x}_i^*} = 0$, otherwise[3]. Clearly, $\delta_{\mathbf{x}_i, \mathbf{x}_i^*} = \delta_{\mathbf{u}_i, \mathbf{u}_i^*} \delta_{\mathbf{g}_i, \mathbf{g}_i^*}$.

---

[2] If $R_i^{st.}(\mathbf{p}, \mathbf{x}_i)$ is some general function of $\mathbf{p}$ and $\mathbf{x}_i$, then $\mathbf{p}$ denotes a set of parameters utilized in a given pricing method.

[3] If the unit status is uniquely set by its output, then instead of $\delta_{\mathbf{x}_i, \mathbf{x}_i^*}$ it is sufficient to use $\delta_{\mathbf{g}_i, \mathbf{g}_i^*}$. Likewise, if the unit output is uniquely set by its status, it suffices to use $\delta_{\mathbf{u}_i, \mathbf{u}_i^*}$ instead of $\delta_{\mathbf{x}_i, \mathbf{x}_i^*}$. The

Let us amend the revenue term in the expression for the profit function by adding some real-valued function $N_i(\mathbf{p},\mathbf{x}_i)$ defined on $R^T \times X_i$. This results in the profit function of the form $\pi_i(\mathbf{p},\mathbf{x}_i) = R_i^{st.}(\mathbf{p},\mathbf{x}_i) + N_i(\mathbf{p},\mathbf{x}_i) - C_i(\mathbf{x}_i) = \pi_i^{st.}(\mathbf{p},\mathbf{x}_i) + N_i(\mathbf{p},\mathbf{x}_i)$. Hence, the generator's profit, obtained when it follows the centralized dispatch outcome, is $\pi_i^*(\mathbf{p}) = \pi_i^{st.}(\mathbf{p},\mathbf{x}_i^*) + N_i(\mathbf{p},\mathbf{x}_i^*)$, while the maximum value of the generator profit function equals $\pi_i^+(\mathbf{p}) = \max_{\mathbf{x}_i \in X_i} \pi_i(\mathbf{p},\mathbf{x}_i)$. In this case, the uplift is expressed as $\pi_i^+(\mathbf{p}) - \pi_i^*(\mathbf{p})$. We impose the following three conditions on $N_i(\mathbf{p},\mathbf{x}_i)$. First, the introduction of $N_i(\mathbf{p},\mathbf{x}_i)$ should not change the maximum generator profit in the decentralized dispatch problem:

$$\pi_i^{st.+}(\mathbf{p}) = \max_{\mathbf{x}_i \in X_i} \pi_i(\mathbf{p},\mathbf{x}_i). \qquad (2)$$

Second, we impose zero uplift condition:

$$\pi_i^{st.+}(\mathbf{p}) = \pi_i^*(\mathbf{p}). \qquad (3)$$

We also require that the new term in the generator revenue function is a rewarding, not penalizing, addition to the standard revenue function $R_i^{st.}(\mathbf{p},\mathbf{x}_i)$:

$$N_i(\mathbf{p},\mathbf{x}_i) \geq 0, \ \forall \mathbf{x}_i \in X_i, \qquad (4)$$

We observe that if the uplift payment is not needed (i.e. $\pi_i^{st.+}(\mathbf{p}) = \pi_i^{st.*}(\mathbf{p})$), then (2) – (4) generally do not yield $N_i(\mathbf{p},\mathbf{x}_i) = 0$, $\forall \mathbf{x}_i \in X_i$, as the profit function may still be amended with no effect on its maximum value and its value at $\mathbf{x}_i^*$. This suggests imposing an additional condition

$$N_i(\mathbf{p},\mathbf{x}_i) = 0, \ \forall \mathbf{x}_i \in X_i, \text{ if } \pi_i^{st.+}(\mathbf{p}) = \pi_i^{st.*}(\mathbf{p}). \ (5)$$

However, this condition can be easily satisfied since given any $N_i'(\mathbf{p},\mathbf{x}_i)$ that satisfies (2) – (4) for $N_i(\mathbf{p},\mathbf{x}_i)$, the function $N_i(\mathbf{p},\mathbf{x}_i) = \theta[\pi_i^{st.+}(\mathbf{p}) - \pi_i^{st.*}(\mathbf{p})]N_i'(\mathbf{p},\mathbf{x}_i)$, with the step-function $\theta(z)$ defined as $\theta(z) = 1$ for $z > 0$ and $\theta(z) = 0$ for $z \leq 0$, satisfies (2) – (5). Therefore, in what follows, we focus on (2) – (4). Since the total profit (including the uplift payment) received by each generator still equals $\pi_i^{st.+}(\mathbf{p})$, it is not affected by $N_i(\mathbf{p},\mathbf{x}_i)$. Thus, the introduction of $N_i(\mathbf{p},\mathbf{x}_i)$ does not address the abovementioned issue of revenue adequacy problem relevant for systems with no fixed load. Also, (2) and (4) entail

$$N(\mathbf{p},\mathbf{x}_i^{st.+}) = 0, \ \forall \mathbf{x}_i^{st.+} \in \arg\max_{\mathbf{x}_i \in X_i} \pi_i^{st.}(\mathbf{p},\mathbf{x}_i), \qquad (6)$$

which gives $\arg\max_{\mathbf{x}_i \in X_i} \pi_i^{st.}(\mathbf{p},\mathbf{x}_i) \subset \arg\max_{\mathbf{x}_i \in X_i} \pi_i(\mathbf{p},\mathbf{x}_i)$. Also, (2) is equivalent to

$$\pi_i(\mathbf{p},\mathbf{x}_i) \leq \pi_i^{st.+}(\mathbf{p}), \ \forall \mathbf{x}_i \in X_i, \qquad (7)$$

$$\pi_i(\mathbf{p},\mathbf{x}_i') = \pi_i^{st.+}(\mathbf{p}), \text{ for some } \mathbf{x}_i' \in X_i. \qquad (8)$$

The obvious choice for $\mathbf{x}_i'$ is $\mathbf{x}_i' = \mathbf{x}_i^*$, which means that the set of (2) and (3) is equivalent to a set of (3) and (7). Another natural choice for $\mathbf{x}_i'$ is $\mathbf{x}_i' = \mathbf{x}_i^{st.+}$, thus (2) is equivalent to a set of (6) and (7). The conditions (2) and (3) imply $\mathbf{x}_i^* \in \arg\max_{\mathbf{x}_i \in X_i} \pi_i(\mathbf{p},\mathbf{x}_i)$, which has the following implication.

---

latter possibility is realized, for example, for a block-loaded unit with the output rate uniquely set by the unit's status.

*Proposition 1* Let $\mathbf{p} \in R^T$ and $R_i^{st.}(\mathbf{p},\mathbf{x}_i) = \mathbf{p}^T\mathbf{g}_i$, $\forall i \in I$. If $\mathbf{x}_i^* \in \arg\max_{\mathbf{x}_i \in X_i} \pi_i(\mathbf{p},\mathbf{x}_i)$, $\forall i \in I$, then $\mathbf{x}^*$ is an optimal point of the following amended centralized dispatch problem

$$\min_{\substack{\mathbf{x}_i \in X_i, \forall i \in I \\ \sum_{i \in I} \mathbf{g}_i = \mathbf{d}}} f_N(\mathbf{x}) \qquad (9)$$

with the objective function $f_N(\mathbf{x}) = \sum_{i \in I}[C_i(\mathbf{x}_i) - N_i(\mathbf{p},\mathbf{x}_i)]$, where $\mathbf{p}$ is treated as the fixed external parameter. Moreover, there is a strong duality between (9) and its dual obtained from the Lagrangian relaxation of the power balance constraint with $\mathbf{p}$ being an optimal value of the dual variable.

Proof. Consider the Lagrangian function $L(\mathbf{q},\mathbf{x}) = \mathbf{q}^T(\mathbf{d} - \sum_{i \in I}\mathbf{g}_i) + f_N(\mathbf{x})$ with a vector of multipliers $\mathbf{q} \in R^T$ and define the dual function $f_N^+(\mathbf{q}) = \min_{\mathbf{x}_i \in X_i, \forall i \in I} L(\mathbf{q},\mathbf{x})$. For $\mathbf{q} = \mathbf{p}$, we have $f_N^+(\mathbf{p}) = \mathbf{p}^T\mathbf{d} - \sum_{i \in I}\max_{\mathbf{x}_i \in X_i} \pi_i(\mathbf{p},\mathbf{x}_i)$. The condition $\mathbf{x}_i^* \in \arg\max_{\mathbf{x}_i \in X_i} \pi_i(\mathbf{p},\mathbf{x}_i)$, $\forall i \in I$, entails $f_N^+(\mathbf{p}) = L(\mathbf{p},\mathbf{x}^*)$. From $L(\mathbf{p},\mathbf{x}^*) = f_N(\mathbf{x}^*)$ we conclude that $f_N^+(\mathbf{p}) = f_N(\mathbf{x}^*)$ - the value of the dual function at $\mathbf{q} = \mathbf{p}$, which is feasible in the dual problem, equals the value of the primal problem objective function at $\mathbf{x} = \mathbf{x}^*$, which is feasible in the primal problem (9). Consequently, we have a strong duality, and ($\mathbf{x}^*$,$\mathbf{p}$) is an optimal primal-dual pair. Proposition is proved.

If conditions of Proposition 1 hold, then $\mathbf{p}$ is a uniform equilibrium price for each generator. Proposition 1 can be straightforwardly generalized to a power market with the price-sensitive consumer bids. In this case, the existence of an equilibrium price does not eliminate the abovementioned budget-balancing problem since (due to the amendments of the consumer cost functions/producer revenue functions) the sum of the consumer payments is at most $\mathbf{p}^T\mathbf{d}$, while the sum of the generator revenues is at least $\mathbf{p}^T\mathbf{d}$.

In the next section, we study the relation between $N_i(\mathbf{p},\mathbf{x}_i)$ satisfying (2) – (4) and the uplift minimization problem.

### III. Utilizing the redundant constraints for the uplift payment reduction

Let us consider some real-valued functions $\rho_i^{l_i}(\mathbf{p},\mathbf{x}_i)$, $l_i \in L_i$, $L_i = \{1,2,...,|L_i|\}$, that are defined on $R^T \times X_i \to R$ and satisfy $\rho_i^{l_i}(\mathbf{p},\mathbf{x}_i) \leq 0$, $\forall \mathbf{x}_i \in X_i$, $\forall \mathbf{p} \in R^T$, $\forall l_i \in L_i$. Introduce a vector function $\rho_i(\mathbf{p},\mathbf{x}_i) = (\rho_i^1(\mathbf{p},\mathbf{x}_i),...,\rho_i^{|L_i|}(\mathbf{p},\mathbf{x}_i))$. Thus, $\rho_i(\mathbf{p},\mathbf{x}_i) \leq 0$, $\forall \mathbf{x}_i \in X_i$, $\forall \mathbf{p} \in R^T$, where we adopt a convention that a vector is non-negative (non-positive) if all of its components are non-negative (non-positive). Clearly, with regard to the centralized dispatch problem (1), the constraints $\rho_i(\mathbf{p},\mathbf{x}_i) \leq 0$, $\forall \mathbf{x}_i \in X_i$, are redundant since they are satisfied on $\Omega$. However, we emphasize that these constraints also hold on a set $\times_{i \in I} X_i$, which contains $\Omega$ as a subset. This means that they belong to a special type of the redundant constraints that hold even if the power balance constraint is removed from the constraint set of (1). (We note that the introduction of the extra copies of generator private equality and/or inequality constraints, which define its private feasible set, also produces this type of redundant constraints.) Among the redundant constraints introduced above there could be the

constraints that hold on $\times_{i \in I} X_i$ as equalities (such as $0 \leq 0$, $u_i^t(u_i^t - 1) \leq 0$) or as strong inequalities (for example, $-1 \leq 0, -\mathbf{g}_i^T \mathbf{g}_i - 1 \leq 0$) - below we show that these kinds of the redundant constraints can be discarded as they do not affect the uplift payment. Consider the optimization problem

$$\max_{\substack{\mathbf{x}_i \in X_i \\ \rho_i(\mathbf{p}, \mathbf{x}_i) \leq 0}} \pi_i(\mathbf{p}, \mathbf{x}_i), \quad (10)$$

which is equivalent to $\max_{\mathbf{x}_i \in X_i} \pi_i(\mathbf{p}, \mathbf{x}_i)$, therefore, $\pi_i^{st.+}(\mathbf{p}) = \max_{\substack{\mathbf{x}_i \in X_i \\ \rho_i(\mathbf{p}, \mathbf{x}_i) \leq 0}} \pi_i(\mathbf{p}, \mathbf{x}_i)$. Let us apply the Lagrangian relaxation procedure to the constraints $\rho_i(\mathbf{p}, \mathbf{x}_i) \leq 0$ and define the Lagrangian function $\pi_i(\mathbf{p}, \mu_i, \mathbf{x}_i) = \pi_i^{st.}(\mathbf{p}, \mathbf{x}_i) - \mu_i^T \rho_i(\mathbf{p}, \mathbf{x}_i)$ with the associated $|L_i|$-dimensional vector of the non-negative multipliers $\mu_i \geq 0$. Clearly, $\pi_i(\mathbf{p}, \mu_i, \mathbf{x}_i)$ is the profit function amended by the non-negative term $-\mu_i^T \rho_i(\mathbf{p}, \mathbf{x}_i)$. Define $\pi_i^+(\mathbf{p}, \mu_i) = \max_{\mathbf{x}_i \in X_i} \pi_i(\mathbf{p}, \mu_i, \mathbf{x}_i)$. Since $\pi_i(\mathbf{p}, \mu_i, \mathbf{x}_i)$ is linear in $\mu_i$, the function $\pi_i^+(\mathbf{p}, \mu_i)$ is convex in $\mu_i$. The problem that is dual to (10) reads:

$$\min_{\mu_i \geq 0} \pi_i^+(\mathbf{p}, \mu_i). \quad (11)$$

*Proposition 2*. There is a strong duality between (10) and (11):

$$\pi_i^{st.+}(\mathbf{p}) = \min_{\mu_i \geq 0} \pi_i^+(\mathbf{p}, \mu_i). \quad (12)$$

Proof. Since $\pi_i^{st.}(\mathbf{p}, \mathbf{x}_i) \leq \pi_i(\mathbf{p}, \mu_i, \mathbf{x}_i)$, $\forall \mathbf{x}_i \in X_i$, $\forall \mu_i \geq 0$, we have $\pi_i^{st.+} = \max_{\mathbf{x}_i \in X_i} \pi_i^{st.}(\mathbf{p}, \mathbf{x}_i) \leq \max_{\mathbf{x}_i \in X_i} \pi_i(\mathbf{p}, \mu_i, \mathbf{x}_i) = \pi_i^+(\mathbf{p}, \mu_i)$, $\forall \mu_i \geq 0$. Hence, $\pi_i^{st.+} \leq \min_{\mu_i \geq 0} \pi_i^+(\mathbf{p}, \mu_i)$. Applying $\min_{\mu_i \geq 0} \pi_i^+(\mathbf{p}, \mu_i) \leq \pi_i^+(\mathbf{p}, 0) = \pi_i^{st.+}$, we obtain $\pi_i^{st.+} = \min_{\mu_i \geq 0} \pi_i^+(\mathbf{p}, \mu_i)$. Proposition is proved.

If $R_i^{st.}(\mathbf{p}, \mathbf{x}_i) = \mathbf{p}^T \mathbf{g}_i$, $\forall i \in I$, with a market price $\mathbf{p}$ obtained using the convex hull pricing method, the same reasoning used to prove Proposition 2 can be applied to the dual problem obtained from (1) by the Lagrangian relaxation of the power balance constraint. Since the redundant constraints under consideration are redundant not only on $\Omega$ but also on $\times_{i \in I} X_i$, the subsequent dualization of these redundant constraints do not affect the duality gap already introduced by dualization of the power balance constraint. Dependence of $\rho_i(\mathbf{p}, \mathbf{x}_i)$ on $\mathbf{p}$ implies that these constraints can be introduced in the dual to the centralized dispatch problem (1) in two different ways. First, $\rho_i(\mathbf{p}, \mathbf{x}_i) \leq 0$ can be added to the constraint set of (1) with some fixed value of $\mathbf{p}$, which is treated as constant in both the primal (1) and the dual problems. Second, the set of constraints $\rho_i(\mathbf{p}, \mathbf{x}_i) \leq 0$ can be introduced after the power balance constraint is relaxed with $\mathbf{p}$ in $\rho_i(\mathbf{p}, \mathbf{x}_i)$ being identified as the multiplier to the power balance constraint in the dual problem (in this case, the dual function to be optimized over $\mathbf{p}$ and $\mu_i$ is generally non-convex in $\mathbf{p}$, but the convexity is restored after optimization over $\mu_i$). In either way, the introduction and subsequent dualization of $\rho_i(\mathbf{p}, \mathbf{x}_i) \leq 0$ (together with the power balance constraint) do not affect the duality gap between (1) and its dual.

In [26] it has been shown that dualization of the affine redundant constraints together with the set of the original constraints of a primal problem results in the same value of the duality gap that emerges from dualization of the original constraints of the primal problem (in our case, this is the power balance constraint), while introduction and dualization of the non-affine redundant constraints may change the value of the dual problem and, hence, affect the duality gap. Although the redundant constraints studied in the present paper are generally non-affine, they do not change the duality gap. The reason is that we deal with a special type of redundant constraints: these constraints hold on $\underset{i \in I}{\times} X_i$, not just on $\Omega$, and Proposition 2 implies that these (possible non-affine) redundant constraints also do not change the duality gap. Thus, the dualization of the redundant constraints, which belong to the specified type, do not affect the value of the dual problem and the duality gap.

Let us define a set $M_i^+(\mathbf{p}) = \arg\min_{\mu_i \in R_{\geq 0}^{L_i}} \pi_i^+(\mathbf{p}, \mu_i)$. Clearly, the set $M_i^+(\mathbf{p})$ generally depends on both the choice of $\rho_i(\mathbf{p}, \mathbf{x}_i)$ and the price $\mathbf{p}$. Since $\pi_i^+(\mathbf{p}, \mu_i)$ is a convex function, the set $M_i^+(\mathbf{p})$ is a closed convex set. Also, (12) entails that $\{0\} \in M_i^+(\mathbf{p})$, which gives that the set $M_i^+(\mathbf{p})$ is nonempty. Due to $\rho_i(\mathbf{p}, \mathbf{x}_i) \leq 0$, $\forall \mathbf{x}_i \in X_i$, and $\mu_i \geq 0$, we have $M_i^+(\mathbf{p}) = \{\mu_i \mid \mu_i \geq 0; \mu_i^T \rho_i(\mathbf{p}, \mathbf{x}_i) \geq \pi_i^{st.}(\mathbf{p}, \mathbf{x}_i) - \pi_i^{st.+}(\mathbf{p}), \forall \mathbf{x}_i \in X_i\}$.

Proposition 2 implies that if $\mu_i^+ \in M_i^+(\mathbf{p})$, then the expression $-\mu_i^{+T} \rho_i(\mathbf{p}, \mathbf{x}_i)$ satisfies conditions (2) and (4) for $N_i(\mathbf{p}, \mathbf{x}_i)$. From (12) we also have a condition on $\mu_i^{+T} \rho_i(\mathbf{p}, \mathbf{x}_i)$ that is equivalent to (6):

$$\mu_i^{+l_i} \rho_i^{l_i}(\mathbf{p}, \mathbf{x}_i^{st.+}) = 0, \ \forall \mathbf{x}_i^{st.+} \in \arg\max_{\mathbf{x}_i \in X_i} \pi_i^{st.}(\mathbf{p}, \mathbf{x}_i), \ \forall l_i \in L_i. \quad (13)$$

Clearly, if for some $l_i \in L_i$ we have $\rho_i^{l_i}(\mathbf{p}, \mathbf{x}_i) = 0$, $\forall \mathbf{x}_i \in X_i$, (i.e. the constraint $\rho_i^{l_i}(\mathbf{p}, \mathbf{x}_i) \leq 0$ is satisfied as equality on $X_i$), then the set of optimal values of $\mu_i^{l_i}$ in (11) is given by $R_{\geq 0}$. The converse is also true: if the set of optimal values of $\mu_i^{l_i}$ is given by $R_{\geq 0}$, then $\rho_i^{l_i}(\mathbf{p}, \mathbf{x}_i) = 0$, $\forall \mathbf{x}_i \in X_i$. Indeed, $\rho_i^{l_i}(\mathbf{p}, \mathbf{x}_i) \neq 0$ for some $\mathbf{x}_i \in X_i$ would imply that the minimum value of $\pi_i^+(\mathbf{p}, \mu_i)$ in (11) can be made arbitrary large, which is impossible for finite $\pi_i^{st.+}(\mathbf{p})$. Likewise, from (13) it follows that if for some $l_i \in L_i$ we have $\rho_i^{l_i}(\mathbf{p}, \mathbf{x}_i) < 0$, $\forall \mathbf{x}_i \in X_i$, (i.e. the constraint $\rho_i^{l_i}(\mathbf{p}, \mathbf{x}_i) \leq 0$ is satisfied as a strong inequality on $X_i$), then the optimal value of $\mu_i^{l_i}$ in (11) is unique and given by $\mu_i^{l_i} = 0$. We note that if either $\rho_i^{l_i}(\mathbf{p}, \mathbf{x}_i) = 0$, $\forall \mathbf{x}_i \in X_i$, or $\rho_i^{l_i}(\mathbf{p}, \mathbf{x}_i) < 0$, $\forall \mathbf{x}_i \in X_i$, then such a constraint makes no contribution to the producer $i$ uplift payment since it does not affect the producer profit.

Clearly, (4) is equivalent to $N_i(\mathbf{p}, \mathbf{x}_i) = -\mu_i^T \rho_i(\mathbf{p}, \mathbf{x}_i)$ with some $\rho_i(\mathbf{p}, \mathbf{x}_i) \leq 0$, $\forall \mathbf{x}_i \in X_i$, and some $\mu_i \geq 0$ (for example, just one constraint $\rho_i(\mathbf{p}, \mathbf{x}_i) = -N_i(\mathbf{p}, \mathbf{x}_i)$ with the multiplier $\mu_i = 1$). The set of equations (2) and (3) for $N_i(\mathbf{p}, \mathbf{x}_i)$ can be transformed to have the form of the optimization problem. Define the producer $i$ uplift payment as $U_i(\mathbf{p}, \mu_i) = \max_{\mathbf{x}_i \in X_i} \pi_i(\mathbf{p}, \mu_i, \mathbf{x}_i) - \pi_i(\mathbf{p}, \mu_i, \mathbf{x}_i^*)$. We have $U_i(\mathbf{p}, \mu_i) \geq 0$, $\forall \mathbf{p} \in R^T$, $\forall \mu_i \geq 0$. Consider the optimization problem

$$\min_{\substack{\mu_i \\ s.t. \\ \mu_i \geq 0, \\ \max_{\mathbf{x}_i \in X_i} \pi_i(\mathbf{p}, \mu_i, \mathbf{x}_i) = \pi_i^{st.+}(\mathbf{p}).}} U_i(\mathbf{p}, \mu_i). \quad (14)$$

Using the definition of $M_i^+(\mathbf{p})$, (14) is expressed as:

$$\min_{\mu_i \in M_i^+(\mathbf{p})} U_i(\mathbf{p},\mu_i) = \pi_i^{st,+}(\mathbf{p}) - \pi_i^{st,*}(\mathbf{p}) + \min_{\mu_i \in M_i^+(\mathbf{p})} \mu_i^T \rho_i(\mathbf{p},\mathbf{x}_i^*). \quad (15)$$

The immediate consequence of (15) is that if $\rho_i(\mathbf{p},\mathbf{x}_i^*) \neq 0$, then $\arg\min_{\mu_i \in M_i^+(\mathbf{p})} U_i(\mathbf{p},\mu_i) \subset \partial M_i^+(\mathbf{p})$, where $\partial M_i^+(\mathbf{p})$ denotes the boundary of $M_i^+(\mathbf{p})$. Thus, the minimum uplift problem for the given price $\mathbf{p}$ and constraint vector function $\rho_i(\mathbf{p},\mathbf{x}_i)$, $\rho_i(\mathbf{p},\mathbf{x}_i^*) \neq 0$, is reduced to the problem of finding a point on the boundary of $M_i^+(\mathbf{p})$ such that the hyperplane containing this point and having the normal vector $-\rho_i(\mathbf{p},\mathbf{x}_i^*)$ supports $M_i^+(\mathbf{p})$, or, equivalently, finding an element of the nonempty closed convex set $M_i^+(\mathbf{p})$ with the largest projection into the direction specified by the vector $-\rho_i(\mathbf{p},\mathbf{x}_i^*)$. We note that for a case of one function $\rho_i(\mathbf{p},\mathbf{x}_i)$ with $\rho_i(\mathbf{p},\mathbf{x}_i^*) \neq 0$, the optimal point of (15) is unique and given by the maximum element of $M_i^+(\mathbf{p})$, i.e. $\mu_i^{+\max}$. We also note that $\min_{\mu_i \in M_i^+(\mathbf{p})} \mu_i^T \rho_i(\mathbf{p},\mathbf{x}_i^*) \leq 0$, which entails $\min_{\mu_i \in M_i^+(\mathbf{p})} U_i(\mathbf{p},\mu_i) \leq \pi_i^{st,+}(\mathbf{p}) - \pi_i^{st,*}(\mathbf{p})$. Thus, addition of $-\mu_i^T \rho_i(\mathbf{p},\mathbf{x}_i)$ to the revenue function results in the lower (or equal) uplift. Clearly, the magnitude of the uplift reduction due to the introduction of $-\mu_i^T \rho_i(\mathbf{p},\mathbf{x}_i)$ in the producer revenue function essentially depends on the choice of $\rho_i(\mathbf{p},\mathbf{x}_i)$, which subsequently specifies the set $M_i^+(\mathbf{p})$. For example, only the redundant constraints that satisfy $\rho_i(\mathbf{p},\mathbf{x}_i) \leq 0$, $\forall \mathbf{x}_i \in X_i$, $\rho_i(\mathbf{p},\mathbf{x}_i^*) \neq 0$, $M_i^+(\mathbf{p}) \neq \{0\}$ may reduce the uplift. The necessary condition for a given vector function $\rho_i(\mathbf{p},\mathbf{x}_i)$ to yield zero uplift payment for the producer $i$ is formulated in the Appendix.

Since $\min_{\mu_i \in M_i^+(\mathbf{p})} U_i(\mathbf{p},\mu_i) \geq 0$, we conclude that $U_i(\mathbf{p},\tilde{\mu}_i) = 0$ for some $\tilde{\mu}_i \in M_i^+(\mathbf{p})$ iff $\min_{\mu_i \in M_i^+(\mathbf{p})} U_i(\mathbf{p},\mu_i) = 0$ and $\tilde{\mu}_i \in \arg\min_{\mu_i \in M_i^+(\mathbf{p})} U_i(\mathbf{p},\mu_i)$. Now we establish a relation between $N_i(\mathbf{p},\mathbf{x}_i)$, which satisfies the conditions (2) - (3), and solutions to (14) with some $\rho_i(\mathbf{p},\mathbf{x}_i) \leq 0$, $\forall \mathbf{x}_i \in X_i$, and $\mu_i \geq 0$.

*Proposition 3.* Let $\rho_i(\mathbf{p},\mathbf{x}_i) \leq 0$, $\forall \mathbf{x}_i \in X_i$, and $\mu_i \geq 0$. The function $N_i(\mathbf{p},\mathbf{x}_i) = -\mu_i^T \rho_i(\mathbf{p},\mathbf{x}_i)$ satisfies the conditions (2) and (3) iff $U_i(\mathbf{p},\mu_i) = 0$ and $\mu_i \in M_i^+(\mathbf{p})$.

Proof. Let the conditions (2) and (3) hold for $N_i(\mathbf{p},\mathbf{x}_i) = -\mu_i^T \rho_i(\mathbf{p},\mathbf{x}_i)$, then (2) implies $\mu_i \in M_i^+(\mathbf{p})$, while (3) yields $U_i(\mathbf{p},\mu_i) = 0$. Likewise, $\mu_i \in M_i^+(\mathbf{p})$ yields (2). Also, both $U_i(\mathbf{p},\mu_i) = 0$ and $\mu_i \in M_i^+(\mathbf{p})$ entail $\min_{\tilde{\mu}_i \in M_i^+(\mathbf{p})} U_i(\mathbf{p},\tilde{\mu}_i) = 0$ and $\mu_i \in \arg\min_{\tilde{\mu}_i \in M_i^+(\mathbf{p})} U_i(\mathbf{p},\tilde{\mu}_i)$, which gives (3) for $N_i(\mathbf{p},\mathbf{x}_i) = -\mu_i^T \rho_i(\mathbf{p},\mathbf{x}_i)$. Proposition is proved.

### IV. Attaining zero uplift payment

Now we formulate different forms of the necessary and sufficient conditions for a given vector function $\rho_i(\mathbf{p},\mathbf{x}_i) \leq 0$, $\forall \mathbf{x}_i \in X_i$, and a multiplier $\mu_i \geq 0$ to produce zero uplift for generator $i$, i.e. for the corresponding $N_i(\mathbf{p},\mathbf{x}_i) = -\mu_i^T \rho_i(\mathbf{p},\mathbf{x}_i)$ to satisfy the conditions (2) and (3). We note that (4) automatically holds for $N_i(\mathbf{p},\mathbf{x}_i) = -\mu_i^T \rho_i(\mathbf{p},\mathbf{x}_i)$.

*Proposition 4.* A function $N_i(\mathbf{p}, \mathbf{x}_i)$ satisfies the conditions (2) - (4) iff $N_i(\mathbf{p}, \mathbf{x}_i) = -\mu_i^T \rho_i(\mathbf{p}, \mathbf{x}_i)$ for with some real-valued vector function $\rho_i(\mathbf{p}, \mathbf{x}_i): X_i \to R^{|L_i|}$ and some multiplier $\mu_i \geq 0$ that satisfy

- $\rho_i(\mathbf{p}, \mathbf{x}_i) \leq 0$, $\forall \mathbf{x}_i \in X_i$; (16)
- $\mu_i^T \rho_i(\mathbf{p}, \mathbf{x}_i^*) = \pi_i^{st,*}(\mathbf{p}) - \pi_i^{st,+}(\mathbf{p})$; (17)
- $\mu_i^T \rho_i(\mathbf{p}, \mathbf{x}_i) \geq \pi_i^{st,\cdot}(\mathbf{p}, \mathbf{x}_i) - \pi_i^{st,+}(\mathbf{p})$, $\forall \mathbf{x}_i \in X_i$. (18)

Proof. Obviously, for a given $N_i(\mathbf{p}, \mathbf{x}_i) = -\mu_i^T \rho_i(\mathbf{p}, \mathbf{x}_i)$ with $\mu_i \geq 0$, the conditions (3) and (4) are equivalent to (17) and (16), respectively. Now we show that (2) holds iff (18) is satisfied, given the validity of (16) and (17). On one hand, if (2) holds with $N_i(\mathbf{p}, \mathbf{x}_i) = -\mu_i^T \rho_i(\mathbf{p}, \mathbf{x}_i)$, then (18) is clearly satisfied. On the other hand, if (18) holds, then $\pi_i^{st,+}(\mathbf{p}) \geq \max_{\mathbf{x}_i \in X_i} [\pi_i^{st,\cdot}(\mathbf{p}, \mathbf{x}_i) - \mu_i^T \rho_i(\mathbf{p}, \mathbf{x}_i)]$ and (17) implies that this weak inequality is satisfied as equality resulting in (2). Proposition is proved.

*Example 1*

For a uninode single-period power market with a market price $p$, let us consider a generator that is offline in a given time interval according to the centralized dispatch solution, i.e. $u_i^* = g_i^* = 0$. Since $\pi_i^{st,*} = 0$ for the offline unit, the generator uplift payment equals $\pi_i^{st,+}(p)$, which is expressed as the generator revenue function amendment of the form $\delta_{x_i, x_i^*} \pi_i^{st,+}(p)$ with $\delta_{x_i, x_i^*} = \delta_{u_i, 0} \delta_{g_i, 0}$. Using $\delta_{u_i, 0} \delta_{g_i, 0} = \delta_{u_i, 0}$, $\forall x_i \in X_i$, and $\delta_{u_i, 0} = 1 - u_i$, the amendment function is expressed as $N_i(p, x_i) = (1 - u_i) \pi_i^{st,+}(p)$, which can be obtained from the redundant constraint $u_i - 1 \leq 0$ with the multiplier $\mu_i = \pi_i^{st,+}(p)$. It can be easily verified that these constraint function and multiplier satisfy (16) – (18).

Combining Propositions 3 and 4 we conclude that the following three statements are equivalent for $\rho_i(\mathbf{p}, \mathbf{x}_i) \leq 0$, $\forall \mathbf{x}_i \in X_i$, and $\tilde{\mu}_i \geq 0$: $\rho_i(\mathbf{p}, \mathbf{x}_i)$ and $\tilde{\mu}_i$ satisfy (17) and (18); $N_i(\mathbf{p}, \mathbf{x}_i) = -\tilde{\mu}_i^T \rho_i(\mathbf{p}, \mathbf{x}_i)$ satisfies (2) - (3); $\min_{\mu_i \in M_i^+(\mathbf{p})} U_i(\mathbf{p}, \mu_i) = 0$ and $\tilde{\mu}_i \in \arg\min_{\mu_i \in M_i^+(\mathbf{p})} U_i(\mathbf{p}, \mu_i)$.

We have proved that (16) - (18) are necessary and sufficient conditions for $N_i(\mathbf{p}, \mathbf{x}_i) = -\mu_i^T \rho_i(\mathbf{p}, \mathbf{x}_i)$ with $\mu_i \geq 0$ to satisfy (2) - (4). Obviously, the choice for $N_i(\mathbf{p}, \mathbf{x}_i)$ (and, therefore, for $\rho_i(\mathbf{p}, \mathbf{x}_i)$ satisfying (16) and $\mu_i \geq 0$) is not unique. Moreover, both the function $\rho_i(\mathbf{p}, \mathbf{x}_i)$ and the multiplier $\mu_i$ that satisfy (16) – (18) may depend on $\mathbf{p}$. (We note that it is always possible to redefine $\rho_i(\mathbf{p}, \mathbf{x}_i)$ so that the resulting $\mu_i$ is independent of $\mathbf{p}$.) Also, it is the scalar product $\mu_i^T \rho_i(\mathbf{p}, \mathbf{x}_i)$, not the individual components $\mu_i^{l_i}$, $\rho_i^{l_i}(\mathbf{p}, \mathbf{x}_i)$, what matters for the uplift calculation.

We note the three important corollaries of Proposition 4. First, if (17) and (18) hold for $\rho_i(\mathbf{p}, \mathbf{x}_i) \leq 0$, $\forall \mathbf{x}_i \in X_i$, and $\mu_i \geq 0$, and for some $l_i \in L_i$ we have $\rho_i^{l_i}(\mathbf{p}, \mathbf{x}_i^*) = 0$, then the same $\rho_i(\mathbf{p}, \mathbf{x}_i)$ with $\mu_i$ modified by setting $\mu_i^{l_i} = 0$ also satisfies (17) and (18). Thus, components of $\rho_i(\mathbf{p}, \mathbf{x}_i)$ that vanish at $\mathbf{x}_i = \mathbf{x}_i^*$ can be excluded from consideration. Second, if (17) and (18) are satisfied by two different pairs $\{\rho_{i|s}(\mathbf{p}, \mathbf{x}_i), \mu_{i|s}\}$: $\rho_{i|s}(\mathbf{p}, \mathbf{x}_i) \leq 0$, $\forall \mathbf{x}_i \in X_i$, $\mu_{i|s} \geq 0$, $s = 1,2$, then so does $\{\rho_{i|3}(\mathbf{p}, \mathbf{x}_i), \mu_{i|3}\}$ with $\mu_{i|3} = 1$ and $\rho_{i|3}(\mathbf{p}, \mathbf{x}_i) = \alpha_1(\mathbf{p}, \mathbf{x}_i) \mu_{i|1}^T \rho_{i|1}(\mathbf{p}, \mathbf{x}_i) + \alpha_2(\mathbf{p}, \mathbf{x}_i) \mu_{i|2}^T \rho_{i|2}(\mathbf{p}, \mathbf{x}_i)$, for any functions $\alpha_s(\mathbf{p}, \mathbf{x}_i)$:

$R^T \times X_i \to R$ that satisfy $\alpha_s(\mathbf{p},\mathbf{x}_i) \geq 0$, $\sum_{s=1,2}\alpha_s(\mathbf{p},\mathbf{x}_i) = 1$, $\forall \mathbf{x}_i \in X_i$. Clearly, $\rho_{i|3}(\mathbf{p},\mathbf{x}_i) \leq 0$, $\forall \mathbf{x}_i \in X_i$. In short, any convex combination of such $\mu_i^T \rho_i(\mathbf{p},\mathbf{x}_i)$ produces another solution for $\rho_i(\mathbf{p},\mathbf{x}_i)$ with $\mu_i = 1$.

Third, if for some set of the redundant constraints $\rho_i'(\mathbf{p},\mathbf{x}_i) \leq 0$, $\forall \mathbf{x}_i \in X_i$, we have found $\mu_i \in M_i^+(\mathbf{p})$, (i.e. the conditions (16) and (18) hold for $\rho_i'(\mathbf{p},\mathbf{x}_i)$ and $\mu_i$), and $\mu_i^T \rho_i'(\mathbf{p},\mathbf{x}_i^*) \neq 0$, then the uplift payment is reduced but is non-zero unless (17) holds. Formally, this can be seen as follows. Let us add to $\rho_i'(\mathbf{p},\mathbf{x}_i)$ an uplift term (i.e. a term proportional to $\delta_{\mathbf{x}_i,\mathbf{x}_i^*}$) that results in (17) while retaining (16) and (18). It is straightforward to verify that (16) - (18) hold for the same $\mu_i \geq 0$ and

$$\rho_i(\mathbf{p},\mathbf{x}_i) = \rho_i'(\mathbf{p},\mathbf{x}_i) + \delta_{\mathbf{x}_i,\mathbf{x}_i^*}[\pi_i^{st,*}(\mathbf{p}) - \pi_i^{st,+}(\mathbf{p}) - \mu_i^T \rho_i'(\mathbf{p},\mathbf{x}_i^*)]\mu_i / \|\mu_i\|. \quad (19)$$

In this case, the function $N_i(\mathbf{p},\mathbf{x}_i)$ is expressed as $N_i(\mathbf{p},\mathbf{x}_i) = -\mu_i^T \rho_i'(\mathbf{p},\mathbf{x}_i) + \delta_{\mathbf{x}_i,\mathbf{x}_i^*}[\pi_i^{st,+}(\mathbf{p}) - \pi_i^{st,*}(\mathbf{p}) + \mu_i^T \rho_i'(\mathbf{p},\mathbf{x}_i^*)]$ with the second term having the form of the uplift payment in the amount of $\pi_i^{st,+}(\mathbf{p}) - \pi_i^{st,*}(\mathbf{p}) + \mu_i^T \rho_i'(\mathbf{p},\mathbf{x}_i^*)$. Clearly, for $\mu_i^T \rho_i'(\mathbf{p},\mathbf{x}_i^*) \neq 0$ we have $0 \leq \mu_i^T \rho_i'(\mathbf{p},\mathbf{x}_i^*) + \pi_i^{st,+}(\mathbf{p}) - \pi_i^{st,*}(\mathbf{p}) < \pi_i^{st,+}(\mathbf{p}) - \pi_i^{st,*}(\mathbf{p})$, which reflects the uplift payment reduction from adding $-\mu_i^T \rho_i'(\mathbf{p},\mathbf{x}_i)$ to the producer revenue function. (Note that addition of $-\mu_i^T \rho_i(\mathbf{p},\mathbf{x}_i)$, with $\rho_i(\mathbf{p},\mathbf{x}_i)$ given by (19), reduces the uplift to zero.) As an illustration, let us observe that $\rho_i'(\mathbf{p},\mathbf{x}_i) = 0$, $\forall \mathbf{x}_i \in X_i$, trivially satisfies (16) and (18) with any $\mu_i > 0$. Utilization of (19) yields $\rho_i(\mathbf{p},\mathbf{x}_i) = \delta_{\mathbf{x}_i,\mathbf{x}_i^*}[\pi_i^{st,*}(\mathbf{p},\mathbf{x}_i^*) - \pi_i^{st,+}(\mathbf{p})]\mu_i / \|\mu_i\|$ resulting in $N_i(\mathbf{p},\mathbf{x}_i) = \delta_{\mathbf{x}_i,\mathbf{x}_i^*}[\pi_i^{st,+}(\mathbf{p}) - \pi_i^{st,*}(\mathbf{p},\mathbf{x}_i^*)]$, which is the revenue function amendment describing the original uplift payment.

Let us denote as $L_i^*$ a subset of $L_i$ with $\rho_i^{l_i}(\mathbf{p},\mathbf{x}_i^*) \neq 0$: $L_i^* = \{l_i \mid l_i \in L_i, \rho_i^{l_i}(\mathbf{p},\mathbf{x}_i^*) \neq 0\}$. This implies that only $\rho_i^{l_i}(\mathbf{p},\mathbf{x}_i)$ with $l_i \in L_i^*$ may contribute to the uplift. Clearly, $L_i^* \subset L_i$. Motivated by Proposition 7, given in Appendix, we have the following statement.

*Proposition 5.* Let $\pi_i^{st,*}(\mathbf{p},\mathbf{x}_i^*) < \pi_i^{st,+}(\mathbf{p})$, $\mu_i \geq 0$, and $\forall l_i \in L_i$ we have $\rho_i^{l_i}(\mathbf{p},\mathbf{x}_i) \leq 0$, $\forall \mathbf{x}_i \in X_i$, then (17) and (18) hold iff

- $L_i^* \neq \emptyset$; (20)
- $\mu_i^{l_i} \leq \inf_{\mathbf{x}_i \in X_{\rho_{l_i} \neq 0}} [\pi_i^{st,*}(\mathbf{p},\mathbf{x}_i) - \pi_i^{st,+}(\mathbf{p}) - \sum_{l_i': l_i' \in L_i, l_i' \neq l_i} \mu_i^{l_i'} \rho_i^{l_i'}(\mathbf{p},\mathbf{x}_i)] / \rho_i^{l_i}(\mathbf{p},\mathbf{x}_i)$, $\forall l_i \in L_i \setminus L_i^*$; (21)
- $\mu_i^{l_i} = \min_{\mathbf{x}_i \in X_{\rho_{l_i} \neq 0}} [\pi_i^{st,*}(\mathbf{p},\mathbf{x}_i) - \pi_i^{st,+}(\mathbf{p}) - \sum_{l_i': l_i' \in L_i, l_i' \neq l_i} \mu_i^{l_i'} \rho_i^{l_i'}(\mathbf{p},\mathbf{x}_i)] / \rho_i^{l_i}(\mathbf{p},\mathbf{x}_i)$, $\forall l_i \in L_i^*$; (22)
- $x_i^* \in \arg\min_{\mathbf{x}_i \in X_{\rho_{l_i} \neq 0}} [\pi_i^{st,*}(\mathbf{p},\mathbf{x}_i) - \pi_i^{st,+}(\mathbf{p}) - \sum_{l_i': l_i' \in L_i, l_i' \neq l_i} \mu_i^{l_i'} \rho_i^{l_i'}(\mathbf{p},\mathbf{x}_i)] / \rho_i^{l_i}(\mathbf{p},\mathbf{x}_i)$, $\forall l_i \in L_i^*$. (23)

Proof. First, we show that, given the assumptions of the proposition, (17) and (18) imply (20) - (23). Indeed, $\pi_i^{st,*}(\mathbf{p},\mathbf{x}_i^*) < \pi_i^{st,+}(\mathbf{p})$ and (17) entail (20). Also, (18) implies $\mu_i^{l_i} \leq \inf_{\mathbf{x}_i \in X_{\rho_{l_i} \neq 0}} [\pi_i^{st,*}(\mathbf{p},\mathbf{x}_i) - \pi_i^{st,+}(\mathbf{p}) - \sum_{l_i': l_i' \in L_i, l_i' \neq l_i} \mu_i^{l_i'} \rho_i^{l_i'}(\mathbf{p},\mathbf{x}_i)] / \rho_i^{l_i}(\mathbf{p},\mathbf{x}_i)$, $\forall l_i \in L_i$, which entails (21) $\forall l_i \in L_i \setminus L_i^*$. For $\forall l_i \in L_i^*$, (17) states that the infimum is attainable at $x_i = x_i^*$ and, therefore, (29) and (30) hold. Second, we prove that under the stated assumptions (20) - (23) entail (17) – (18). Clearly, (21) and (22) give

$$\mu_i^{l_i} \le \inf_{\mathbf{x}_i \in X_{\rho_i^{l_i} \ne 0}} [\pi_i^{st.}(\mathbf{p},\mathbf{x}_i) - \pi_i^{st.+}(\mathbf{p}) - \sum_{l_i': l_i' \in L_i, l_i' \ne l_i} \mu_i^{l_i'} \rho_i^{l_i'}(\mathbf{p},\mathbf{x}_i)] / \rho_i^{l_i}(\mathbf{p},\mathbf{x}_i), \quad \forall l_i \in L_i,$$

which implies

$$\mu_i^T \rho_i(\mathbf{p},\mathbf{x}_i) \ge \pi_i^{st.}(\mathbf{p},\mathbf{x}_i) - \pi_i^{st.+}(\mathbf{p}), \quad \forall \mathbf{x}_i \in \bigcup_{l_i \in L_i} X_{\rho_i^{l_i} \ne 0}.$$

Since $\mu_i^T \rho_i(\mathbf{p},\mathbf{x}_i) = 0$, $\forall \mathbf{x}_i \in X_i \setminus \bigcup_{l_i \in L_i} X_{\rho_i^{l_i} \ne 0}$, we conclude that (18) holds. Due to (20) and (22), for some value of $l_i$ we have

$$\mu_i^{l_i} = \min_{\mathbf{x}_i \in X_{\rho_i^{l_i} \ne 0}} [\pi_i^{st.}(\mathbf{p},\mathbf{x}_i) - \pi_i^{st.+}(\mathbf{p}) - \sum_{l_i': l_i' \in L_i, l_i' \ne l_i} \mu_i^{l_i'} \rho_i^{l_i'}(\mathbf{p},\mathbf{x}_i)] / \rho_i^{l_i}(\mathbf{p},\mathbf{x}_i),$$

$$x_i^* \in \arg\min_{\mathbf{x}_i \in X_{\rho_i^{l_i} \ne 0}} [\pi_i^{st.}(\mathbf{p},\mathbf{x}_i) - \pi_i^{st.+}(\mathbf{p}) - \sum_{l_i': l_i' \in L_i, l_i' \ne l_i} \mu_i^{l_i'} \rho_i^{l_i'}(\mathbf{p},\mathbf{x}_i)] / \rho_i^{l_i}(\mathbf{p},\mathbf{x}_i),$$

which entails (17). Proposition is proved.

We note that to derive (17) from (20), (22), and (23) we needed validity of (22) and (23) for just one $l_i \in L_i^*$. This is because the set (21) – (23) is equivalent to the set of the following statements: (21) is valid for $\forall l_i \in L_i$, (22) and (23) hold for some $l_i \in L_i^*$.

Thus, Propositions 4, 5 give the necessary and sufficient conditions for a given $\rho_i(\mathbf{p},\mathbf{x}_i) \le 0$, $\forall \mathbf{x}_i \in X_i$, and $\mu_i \ge 0$ to yield zero uplift payment for the producer $i$.

We also note that the different formulations of the redundant constraints generally result in the non-equivalent amendments of the generator revenue function. For example, the constraint set $\rho_i^1(\mathbf{p},\mathbf{x}_i) \le 0$, $\rho_i^2(\mathbf{p},\mathbf{x}_i) \le 0$ is equivalent to $\max[\rho_i^1(\mathbf{p},\mathbf{x}_i), \rho_i^2(\mathbf{p},\mathbf{x}_i)] \le 0$. However, the function $-\mu \max[\rho_i^1(\mathbf{p},\mathbf{x}_i), \rho_i^2(\mathbf{p},\mathbf{x}_i)]$ with some multiplier $\mu \ge 0$ generally cannot be expressed as $-\mu^1 \rho_i^1(\mathbf{p},\mathbf{x}_i) - \mu^2 \rho_i^2(\mathbf{p},\mathbf{x}_i)$ with some multipliers $\mu^1 \ge 0$, $\mu^2 \ge 0$. In certain cases, the redundant constraint $\max[\rho_i^1(\mathbf{p},\mathbf{x}_i), \rho_i^2(\mathbf{p},\mathbf{x}_i)] \le 0$ satisfies (16) – (18) with some $\mu \ge 0$, while the constraints $\rho_i^1(\mathbf{p},\mathbf{x}_i) \le 0$, $\rho_i^2(\mathbf{p},\mathbf{x}_i) \le 0$ do not satisfy (16) – (18) for any $\mu^1 \ge 0$, $\mu^2 \ge 0$. This is because the transition from a set of the constraints $\rho_i^1(\mathbf{p},\mathbf{x}_i) \le 0$, $\rho_i^2(\mathbf{p},\mathbf{x}_i) \le 0$ to the equivalent constraint $\max[\rho_i^1(\mathbf{p},\mathbf{x}_i), \rho_i^2(\mathbf{p},\mathbf{x}_i)] \le 0$ is a nonlinear operation, while the considered amendment functions are linear in the redundant constraints.

Clearly, just one appropriate redundant constraint is sufficient to obtain zero uplift payment for the generator $i$: for example, $\rho_i(\mathbf{p},\mathbf{x}_i) = -N_i(\mathbf{p},\mathbf{x}_i)$ with the associated multiplier equal 1. Likewise, given a vector function $\rho_i(\mathbf{p},\mathbf{x}_i)$ and the associated vector of the non-negative multipliers $\mu_i \ge 0$ that satisfy (16) – (18), just one redundant constraint of the form $-\mu_i^T \rho_i(\mathbf{p},\mathbf{x}_i) \le 0$ yields zero uplift payment for the producer.

Moreover, in case of the uniform market price $\mathbf{p} \in R^T$ (which implies $R_i^{st.}(\mathbf{p},\mathbf{x}_i) = \mathbf{p}^T \mathbf{g}_i$, $\forall i \in I$), it suffices to introduce only one redundant constraint in (1) to obtain zero total uplift payment. Indeed, for a given set of $N_i(\mathbf{p},\mathbf{x}_i)$, $i \in I$, each satisfying (2) – (4), let us introduce the redundant constraint $-\sum_{i \in I} N_i(\mathbf{p},\mathbf{x}_i) \le 0$ in the centralized dispatch optimization problem (1). Clearly, this constraint affects neither the value of (1) nor the set of its optimal points. Let us apply the Lagrangian relaxation procedure to both this constraint and the power balance constraint with the multipliers $\nu \geq 0$ and $\mathbf{q} \in R^T$, respectively. This yields $L(\mathbf{q},\nu,\mathbf{x}) = \mathbf{q}^T(\mathbf{d} - \sum_{i \in I} \mathbf{g}_i) + \sum_{i \in I} C_i(\mathbf{x}_i) - \nu \sum_{i \in I} N_i(\mathbf{p},\mathbf{x}_i)$, which entails

$$\min_{\mathbf{x}_i \in X_i, \forall i \in I} L(\mathbf{q},\nu,\mathbf{x}) = \mathbf{q}^T\mathbf{d} - \sum_{i \in I} \max_{\mathbf{x}_i \in X_i}[\pi_i^{st\cdot}(\mathbf{q},\mathbf{x}_i) + \nu N_i(\mathbf{p},\mathbf{x}_i)] \leq f^* - \nu \sum_{i \in I} N_i(\mathbf{p},\mathbf{x}_i^*).$$

For a given $\mathbf{q} \in R^T$ and $\nu \geq 0$, the total uplift payment is

$$\sum_{i \in I}[\max_{\mathbf{x}_i \in X_i}[\pi_i^{st\cdot}(\mathbf{q},\mathbf{x}_i) + \nu N_i(\mathbf{p},\mathbf{x}_i)] - \pi_i^{st\cdot}(\mathbf{p},\mathbf{x}_i^*) - \nu N_i(\mathbf{p},\mathbf{x}_i^*)] = f^* - \min_{\mathbf{x}_i \in X_i, \forall i \in I} L(\mathbf{p},\nu,\mathbf{x}) - \nu \sum_{i \in I} N_i(\mathbf{p},\mathbf{x}_i^*).$$

For $\mathbf{q}=\mathbf{p}$ and $\nu=1$, (2) implies $\min_{\mathbf{x}_i \in X_i, \forall i \in I} L(\mathbf{p},1,\mathbf{x}) = \mathbf{p}^T\mathbf{d} - \sum_{i \in I}\pi_i^{st\cdot+}(\mathbf{p}) = \min_{\mathbf{x}_i \in X_i, \forall i \in I} L(\mathbf{p},0,\mathbf{x})$.

Consequently, the total uplift payment equals

$$\sum_{i \in I}[\pi_i^{st\cdot+}(\mathbf{p}) - \pi_i^{st\cdot}(\mathbf{p},\mathbf{x}_i^*) - N_i(\mathbf{p},\mathbf{x}_i^*)] = f^* - \min_{\mathbf{x}_i \in X_i, \forall i \in I} L(\mathbf{p},0,\mathbf{x}) - \sum_{i \in I} N_i(\mathbf{p},\mathbf{x}_i^*),$$

which is zero due to (3).

If the uniform market price $\mathbf{p}$ is obtained using the convex hull pricing method, then $f^* - \min_{\mathbf{x}_i \in X_i, \forall i \in I} L(\mathbf{p},0,\mathbf{x})$ is the original duality gap, but the total uplift payment no longer equals the duality gap and can be reduced to zero by utilizing the proper functions $N_i(\mathbf{p},\mathbf{x}_i)$, $i \in I$. We also note the following relation to the dual problem. From (4) it follows that $\min_{\nu \geq 0}\max_{\mathbf{x}_i \in X_i}[\pi_i^{st\cdot}(\mathbf{q},\mathbf{x}_i) + \nu N_i(\mathbf{p},\mathbf{x}_i)] = \max_{\mathbf{x}_i \in X_i}\pi_i^{st\cdot}(\mathbf{q},\mathbf{x}_i)$, $\forall \mathbf{q} \in R^T$, $\forall i \in I$. Therefore,

$$L(\mathbf{q},0,\mathbf{x}) \leq \max_{\nu \geq 0}\min_{\mathbf{x}_i \in X_i, \forall i \in I} L(\mathbf{q},\nu,\mathbf{x}) = \mathbf{q}^T\mathbf{d} - \min_{\nu \geq 0}\sum_{i \in I}\max_{\mathbf{x}_i \in X_i}[\pi_i^{st\cdot}(\mathbf{q},\mathbf{x}_i) + \nu N_i(\mathbf{p},\mathbf{x}_i)] \leq \mathbf{q}^T\mathbf{d}$$
$$- \sum_{i \in I}\min_{\nu \geq 0}\max_{\mathbf{x}_i \in X_i}[\pi_i^{st\cdot}(\mathbf{q},\mathbf{x}_i) + \nu N_i(\mathbf{p},\mathbf{x}_i)] = L(\mathbf{q},0,\mathbf{x}).$$

Thus, $\max_{\mathbf{q} \in R^T}\max_{\nu \geq 0}\min_{\mathbf{x}_i \in X_i, \forall i \in I} L(\mathbf{q},\nu,\mathbf{x}) = \max_{\mathbf{q} \in R^T}\min_{\mathbf{x}_i \in X_i, \forall i \in I} L(\mathbf{q},0,\mathbf{x})$, and neither the value of the dual problem that is obtained from (1) by dualizing the power balance constraint nor the set of the optimal dual variables $\mathbf{q}$ is affected by the introduction of the redundant constraint. Consequently, if the market price is set by the convex hull pricing method, then the redundant constraint $-\sum_{i \in I} N_i(\mathbf{p},\mathbf{x}_i) \leq 0$ yields the same optimal set of the market prices and maximum values of the producer profit $\pi_i^{st\cdot}(\mathbf{p},\mathbf{x}_i)$, $\forall i \in I$, and results in zero total uplift payment.

### V. General form expression for $N_i(\mathbf{p},\mathbf{x}_i)$

We note that for a given $\mathbf{p}$ the function $N_i(\mathbf{p},\mathbf{x}_i)$ that satisfies (2) and (4) takes values between zero and the non-negative function $\pi_i^{st\cdot+}(\mathbf{p}) - \pi_i^{st\cdot}(\mathbf{p},\mathbf{x}_i)$, while (3) fixes its value at $\mathbf{x}_i^*$ by $N_i(\mathbf{p},\mathbf{x}_i^*) = \pi_i^{st\cdot+}(\mathbf{p}) - \pi_i^{st\cdot*}(\mathbf{p})$. Motivated by this observation we have the following general form expression for $N_i(\mathbf{p},\mathbf{x}_i) = -\mu_i^T \rho_i(\mathbf{p},\mathbf{x}_i)$.

*Proposition 6.* The conditions (2) - (4) hold for $N_i(\mathbf{p},\mathbf{x}_i)$ iff $N_i(\mathbf{p},\mathbf{x}_i)$ satisfies

$$N_i(\mathbf{p},\mathbf{x}_i) = \min[\pi_i^{st\cdot+}(\mathbf{p}) - \pi_i^{st\cdot}(\mathbf{p},\mathbf{x}_i); \delta_{\mathbf{x}_i,\mathbf{x}_i^*}[\pi_i^{st\cdot+}(\mathbf{p}) - \pi_i^{st\cdot*}(\mathbf{p})] + \gamma_i(\mathbf{p},\mathbf{x}_i)], \forall \mathbf{x}_i \in X_i, \quad (24)$$

with some non-negative real-valued function $\gamma_i(\mathbf{p},\mathbf{x}_i): R^T \times X_i \to R$, $\gamma_i(\mathbf{p},\mathbf{x}_i) \geq 0$, $\forall \mathbf{x}_i \in X_i$.

Proof. It is straightforward to verify that (24) satisfies (2) – (4). Now we show that (2) – (4) entail (24) with some function $\gamma_i(\mathbf{p},\mathbf{x}_i): R^T \times X_i \to R$, $\gamma_i(\mathbf{p},\mathbf{x}_i) \geq 0$, $\forall \mathbf{x}_i \in X_i$. Consider a function $\gamma_i(\mathbf{p},\mathbf{x}_i) = N_i(\mathbf{p},\mathbf{x}_i) - \delta_{\mathbf{x}_i,\mathbf{x}_i^*}[\pi_i^{st.+}(\mathbf{p}) - \pi_i^{st.*}(\mathbf{p})]$, $\forall \mathbf{x}_i \in X_i$. Since $N_i(\mathbf{p},\mathbf{x}_i)$ satisfies (3) and (4), we have $\gamma_i(\mathbf{p},\mathbf{x}_i) \geq 0$, $\forall \mathbf{x}_i \in X_i$. From (2) we have $N_i(\mathbf{p},\mathbf{x}_i) \leq \pi_i^{st.+}(\mathbf{p}) - \pi_i^{st.}(\mathbf{p},\mathbf{x}_i)$, $\forall \mathbf{x}_i \in X_i$, which gives $N_i(\mathbf{p},\mathbf{x}_i) = \min[\pi_i^{st.+}(\mathbf{p}) - \pi_i^{st.}(\mathbf{p},\mathbf{x}_i); N_i(\mathbf{p},\mathbf{x}_i)]$, $\forall \mathbf{x}_i \in X_i$. Using $N_i(\mathbf{p},\mathbf{x}_i) = \gamma_i(\mathbf{p},\mathbf{x}_i) + \delta_{\mathbf{x}_i,\mathbf{x}_i^*}[\pi_i^{st.+}(\mathbf{p}) - \pi_i^{st.*}(\mathbf{p})]$, we obtain (24). Proposition is proved.

From (24) we have the following general expression for the amended profit function

$$\pi_i(\mathbf{p},\mathbf{x}_i) = \min[\pi_i^{st.+}(\mathbf{p}); \pi_i^{st.}(\mathbf{p},\mathbf{x}_i) + \delta_{\mathbf{x}_i,\mathbf{x}_i^*}[\pi_i^{st.+}(\mathbf{p}) - \pi_i^{st.*}(\mathbf{p})] + \gamma_i(\mathbf{p},\mathbf{x}_i)] \qquad (25)$$

parameterized by the function $\gamma_i(\mathbf{p},\mathbf{x}_i)$, which is non-negative on $X_i$. From (25) it follows that the resulting amended profit function on $X_i$ majorizes the function $\pi_i^{st.}(\mathbf{p},\mathbf{x}_i) + \delta_{\mathbf{x}_i,\mathbf{x}_i^*}[\pi_i^{st.+}(\mathbf{p}) - \pi_i^{st.*}(\mathbf{p})]$ and is bounded by $\pi_i^{st.+}(\mathbf{p})$. The converse is also true: any function on $X_i$ that majorizes $\pi_i^{st.}(\mathbf{p},\mathbf{x}_i) + \delta_{\mathbf{x}_i,\mathbf{x}_i^*}[\pi_i^{st.+}(\mathbf{p}) - \pi_i^{st.*}(\mathbf{p})]$ and is bounded by $\pi_i^{st.+}(\mathbf{p})$ satisfies (25) with some $\gamma_i(\mathbf{p},\mathbf{x}_i) \geq 0$, $\forall \mathbf{x}_i \in X_i$. It is worth mentioning that (24) implies (6). Now we study some forms of $N_i(\mathbf{p},\mathbf{x}_i)$ generated by various choices of the parameter function $\gamma(\mathbf{p},\mathbf{x}_i)$.

*Example 2*

Setting $\gamma_i(\mathbf{p},\mathbf{x}_i) = 0$, $\forall \mathbf{x}_i \in X_i$, in (24) gives $N_i(\mathbf{p},\mathbf{x}_i) = \delta_{\mathbf{x}_i,\mathbf{x}_i^*}[\pi_i^{st.+}(\mathbf{p}) - \pi_i^{st.}(\mathbf{p},\mathbf{x}_i^*)]$, which is the original uplift payment of $\pi_i^{st.+}(\mathbf{p}) - \pi_i^{st.}(\mathbf{p},\mathbf{x}_i^*)$ recast in the form of the revenue function amendment.

This expression for $N_i(\mathbf{p},\mathbf{x}_i)$ can be easily obtained from the redundant constraint $\rho_i(\mathbf{p},\mathbf{x}_i) = -\delta_{\mathbf{x}_i,\mathbf{x}_i^*}$, which satisfies $\rho_i(\mathbf{p},\mathbf{x}_i) \leq 0$, $\forall \mathbf{x}_i \in X_i$. Application of Proposition 7 gives $\mu_i^{+\max} = \pi_i^{st.+}(\mathbf{p}) - \pi_i^{st.}(\mathbf{p},\mathbf{x}_i^*)$. Clearly, these $\rho_i(\mathbf{p},\mathbf{x}_i)$ and $\mu_i^{+\max}$ satisfy (16) – (18). Hence, $N_i(\mathbf{p},\mathbf{x}_i) = \delta_{\mathbf{x}_i,\mathbf{x}_i^*}[\pi_i^{st.+}(\mathbf{p}) - \pi_i^{st.}(\mathbf{p},\mathbf{x}_i^*)]$ satisfies (2) – (4).

*Example 3*

Setting $\gamma_i(\mathbf{p},\mathbf{x}_i) = (1 - \delta_{\mathbf{x}_i,\mathbf{x}_i^*})[\pi_i^{st.+}(\mathbf{p}) - \pi_i^{st.}(\mathbf{p},\mathbf{x}_i)]$ in (24) results in $N_i(\mathbf{p},\mathbf{x}_i) = \pi_i^{st.+}(\mathbf{p}) - \pi_i^{st.}(\mathbf{p},\mathbf{x}_i)$, which implies $\pi_i(\mathbf{p},\mathbf{x}_i) = \pi_i^{st.+}(\mathbf{p})$, $\forall \mathbf{x}_i \in X_i$. In this case, the profit function becomes a constant independent of $\mathbf{x}_i$ and the producer is indifferent to its output volume. (If unacceptable, such a solution can be easily excluded by adding a condition that $N_i(\mathbf{p},\mathbf{x}_i) < \pi_i^{st.+}(\mathbf{p}) - \pi_i^{st.}(\mathbf{p},\mathbf{x}_i)$ for some $\mathbf{x}_i \in X_i$.)

*Example 4*

Let us utilize the freedom to choose an arbitrary (non-negative on $X_i$) function $\gamma_i(\mathbf{p},\mathbf{x}_i)$ in (25) to smooth the possible discontinuity of $\pi_i(\mathbf{p},\mathbf{x}_i)$, as a function of $\mathbf{x}_i$, introduced by the discontinuous term $\delta_{\mathbf{x}_i,\mathbf{x}_i^*}[\pi_i^{st.+}(\mathbf{p}) - \pi_i^{st.*}(\mathbf{p})]$. Consider the extended value function $\bar{\pi}_i^{st.}(\mathbf{p},\mathbf{x}_i)$ on $Conv\{X_i\}$ defined as $\bar{\pi}_i^{st.}(\mathbf{p},\mathbf{x}_i) = \pi_i^{st.}(\mathbf{p},\mathbf{x}_i)$ for $\mathbf{x}_i \in X_i$ and $\bar{\pi}_i^{st.}(\mathbf{p},\mathbf{x}_i) = -\infty$ for $\mathbf{x}_i \in Conv\{X_i\} \setminus X_i$. Let us choose

$$\gamma_i(\mathbf{p},\mathbf{x}_i) = conc_{Conv\{X_i\}}\{\overline{\pi}_i^{st.}(\mathbf{p},\mathbf{x}_i) + \delta_{\mathbf{x}_i,\mathbf{x}_i^*}[\pi_i^{st.+}(\mathbf{p}) - \pi_i^{st.*}(\mathbf{p})]\} - \pi_i^{st.}(\mathbf{p},\mathbf{x}_i)$$
$$- \delta_{\mathbf{x}_i,\mathbf{x}_i^*}[\pi_i^{st.+}(\mathbf{p}) - \pi_i^{st.*}(\mathbf{p})],$$

where $conc_{Conv\{X_i\}}$ denotes the concave hull of a function on the set $Conv\{X_i\}$. Clearly, $\gamma_i(\mathbf{p},\mathbf{x}_i) \geq 0$, $\forall \mathbf{x}_i \in X_i$. Therefore, (25) entails

$$\pi_i(\mathbf{p},\mathbf{x}_i) = \min[\pi_i^{st.+}(\mathbf{p}); conc_{Conv\{X_i\}}\{\overline{\pi}_i^{st.}(\mathbf{p},\mathbf{x}_i) + \delta_{\mathbf{x}_i,\mathbf{x}_i^*}[\pi_i^{st.+}(\mathbf{p}) - \pi_i^{st.*}(\mathbf{p})]\}].$$

From $\pi_i^{st.}(\mathbf{p},\mathbf{x}_i) + \delta_{\mathbf{x}_i,\mathbf{x}_i^*}[\pi_i^{st.+}(\mathbf{p}) - \pi_i^{st.*}(\mathbf{p})] \leq \pi_i^{st.+}(\mathbf{p})$, $\forall \mathbf{x}_i \in X_i$, and $\max_{\mathbf{x}_i \in X_i}\{\pi_i^{st.}(\mathbf{p},\mathbf{x}_i) + \delta_{\mathbf{x}_i,\mathbf{x}_i^*}[\pi_i^{st.+}(\mathbf{p}) - \pi_i^{st.*}(\mathbf{p})]\} = \max_{\mathbf{x}_i \in X_i} conc_{Conv\{X_i\}}\{\overline{\pi}_i^{st.}(\mathbf{p},\mathbf{x}_i) + \delta_{\mathbf{x}_i,\mathbf{x}_i^*}[\pi_i^{st.+}(\mathbf{p}) - \pi_i^{st.*}(\mathbf{p})]\}$

we obtain $conc_{Conv\{X_i\}}\{\overline{\pi}_i^{st.}(\mathbf{p},\mathbf{x}_i) + \delta_{\mathbf{x}_i,\mathbf{x}_i^*}[\pi_i^{st.+}(\mathbf{p}) - \pi_i^{st.*}(\mathbf{p})]\} \leq \pi_i^{st.+}(\mathbf{p})$, $\forall \mathbf{x}_i \in X_i$. This results in

$$\pi_i(\mathbf{p},\mathbf{x}_i) = conc_{Conv\{X_i\}}\{\overline{\pi}_i^{st.}(\mathbf{p},\mathbf{x}_i) + \delta_{\mathbf{x}_i,\mathbf{x}_i^*}[\pi_i^{st.+}(\mathbf{p}) - \pi_i^{st.*}(\mathbf{p})]\}, \forall \mathbf{x}_i \in X_i.$$

Clearly, since the function $\pi_i(\mathbf{p},\mathbf{x}_i)$ is concave on $Conv\{X_i\}$, it is continuous in $\mathbf{x}_i$ in the interior of $Conv\{X_i\}$. For $R_i^{st.}(\mathbf{p},\mathbf{x}_i) = \mathbf{p}^T\mathbf{g}_i$, it is straightforward to verify that $\forall \mathbf{x}_i \in Conv\{X_i\}$ we have

$$conc_{Conv\{X_i\}}\{\overline{\pi}_i^{st.}(\mathbf{p},\mathbf{x}_i) + \delta_{\mathbf{x}_i,\mathbf{x}_i^*}[\pi_i^{st.+}(\mathbf{p}) - \pi_i^{st.*}(\mathbf{p})]\} = \mathbf{p}^T\mathbf{g}_i - conv_{Conv\{X_i\}}\{\overline{C}_i(\mathbf{x}_i)$$
$$- \delta_{\mathbf{x}_i,\mathbf{x}_i^*}[\pi_i^{st.+}(\mathbf{p}) - \pi_i^{st.*}(\mathbf{p})]\},$$

where $\overline{C}_i(\mathbf{x}_i)$ is the extended value function defined on $Conv\{X_i\}$ as $\overline{C}_i(\mathbf{x}_i) = C_i(\mathbf{x}_i)$ for $\mathbf{x}_i \in X_i$ and $\overline{C}_i(\mathbf{x}_i) = +\infty$ for $\mathbf{x}_i \in Conv\{X_i\} \setminus X_i$, and $conv_{Conv\{X_i\}}$ denotes the convex hull of a function on the set $Conv\{X_i\}$. Therefore, in case of $R_i^{st.}(\mathbf{p},\mathbf{x}_i) = \mathbf{p}^T\mathbf{g}_i$, we have

$$\pi_i(\mathbf{p},\mathbf{x}_i) = \mathbf{p}^T\mathbf{g}_i - conv_{Conv\{X_i\}}\{\overline{C}_i(\mathbf{x}_i) - \delta_{\mathbf{x}_i,\mathbf{x}_i^*}[\pi_i^{st.+}(\mathbf{p}) - \pi_i^{st.*}(\mathbf{p})]\}, \forall \mathbf{x}_i \in X_i, \quad (26)$$

which gives the following expression for $N_i(\mathbf{p},\mathbf{x}_i)$:

$$N_i(\mathbf{p},\mathbf{x}_i) = C_i(\mathbf{x}_i) - conv_{Conv\{X_i\}}\{\overline{C}_i(\mathbf{x}_i) - \delta_{\mathbf{x}_i,\mathbf{x}_i^*}[\pi_i^{st.+}(\mathbf{p}) - \pi_i^{st.*}(\mathbf{p})]\}, \quad \forall \mathbf{x}_i \in X_i. \quad (27)$$

Clearly, (27) can be realized using one redundant constraint $conv_{Conv\{X_i\}}\{\overline{C}_i(\mathbf{x}_i) - \delta_{\mathbf{x}_i,\mathbf{x}_i^*}[\pi_i^{st.+}(\mathbf{p}) - \pi_i^{st.*}(\mathbf{p})]\} - C_i(\mathbf{x}_i) \leq 0$, $\forall \mathbf{x}_i \in X_i$, with the associated multiplier $\mu_i = 1$.

## VI. Application to power markets with marginal pricing

In the case of a power market with marginal pricing, $\mathbf{p}$ is identified as the marginal price faced by the generator $i$ (the system marginal price or the locational marginal price at the generator node), and $R_i^{st.}(\mathbf{p},\mathbf{x}_i) = \mathbf{p}^T\mathbf{g}_i$. Let us assume that for each fixed value of $\mathbf{u}_i$ both the generator private feasible set and the cost function are convex. For the market price $\mathbf{p}$ and a given fixed vector of statuses $\mathbf{u}_i' \in \{0,1\}^T$, let us introduce the maximum value of the generator $i$ profit $\pi_i^{st.\max}(\mathbf{p},\mathbf{u}_i')$:

$$\pi_i^{st.\max}(\mathbf{p},\mathbf{u}_i') = \max_{\substack{\mathbf{x}_i \\ s.t. \\ \mathbf{x}_i \in X_i, \\ \mathbf{u}_i = \mathbf{u}_i'.}} \pi_i^{st.}(\mathbf{p},\mathbf{x}_i). \quad (28)$$

The marginal pricing method entails

$$\pi_i^{st.*}(\mathbf{p}) = \pi_i^{st.\max}(\mathbf{p},\mathbf{u}_i^*). \quad (29)$$

As it was mentioned above, the choice for $N_i(\mathbf{p},\mathbf{x}_i)$ satisfying (2) – (4) is not unique. Below we provide two different expressions for $N_i(\mathbf{p},\mathbf{x}_i)$, which originate from the two different sets of the redundant constraints.

For each time instance $t \in \{1,...,T\}$, we have the redundant inequality constraints on the status binary variable: $u_i^t \geq 0$ and $1 - u_i^t \geq 0$. Let us choose the constraint $1 - u_i^t \geq 0$ if $u_i^{*t} = 0$, and the constraint $u_i^t \geq 0$ if $u_i^{*t} = 1$. Such a choice can be written as $(u_i^t)^{u_i^{*t}}(1-u_i^t)^{(1-u_i^{*t})} \geq 0$. Consider the redundant constraint that results from the product of these constraints for all the time instances: $\rho_i(\mathbf{u}_i) \leq 0$, $\forall \mathbf{x}_i \in X_i$, with $\rho_i(\mathbf{u}_i) = -\prod_{t \in \{1,...T\}}(u_i^t)^{u_i^{*t}}(1-u_i^t)^{(1-u_i^{*t})}$. It is straightforward to check that $\rho_i(\mathbf{u}_i) = 0$ if $\mathbf{u}_i \neq \mathbf{u}_i^*$, and $\rho_i(\mathbf{u}_i) = -1$ if $\mathbf{u}_i = \mathbf{u}_i^*$. Therefore, $\rho_i(\mathbf{u}_i) = -\delta_{\mathbf{u}_i,\mathbf{u}_i^*}$. Proposition 7 gives the value of the associated multiplier $\mu_i^{+\max} = \pi_i^{st.+}(\mathbf{p}) - \pi_i^{st.\max}(\mathbf{p},\mathbf{u}_i^*)$. Using (29), we arrive at $\mu_i^{+\max} = \pi_i^{st.+}(\mathbf{p}) - \pi_i^{st.*}(\mathbf{p})$. Consequently, the redundant constraint $\rho_i(\mathbf{u}_i) = -\delta_{\mathbf{u}_i,\mathbf{u}_i^*}$ and the associated multiplier $\mu_i^{+\max}$ satisfy (16) – (18). Thus, the conditions (2) – (4) hold for $N_i(\mathbf{p},\mathbf{x}_i) = \delta_{\mathbf{u}_i,\mathbf{u}_i^*}[\pi_i^{st.+}(\mathbf{p}) - \pi_i^{st.*}(\mathbf{p})]$. Indeed, $N_i(\mathbf{p},\mathbf{x}_i)$ is non-negative, has the right value at $\mathbf{x}_i = \mathbf{x}_i^*$ and has no effect on the maximum value of the profit function:

$$\max_{\mathbf{x}_i \in X_i}[\pi_i^{st.}(\mathbf{p},\mathbf{x}_i) + N_i(\mathbf{p},\mathbf{x}_i)] = \max\{\max_{\substack{\mathbf{x}_i \\ s.t. \\ \mathbf{x}_i \in X_i, \\ \mathbf{u}_i \neq \mathbf{u}_i^*}} \pi_i^{st.}(\mathbf{p},\mathbf{x}_i); \pi_i^{st.+}(\mathbf{p}) - \pi_i^{st.*}(\mathbf{p}) + \max_{\substack{\mathbf{x}_i \\ s.t. \\ \mathbf{x}_i \in X_i, \\ \mathbf{u}_i = \mathbf{u}_i^*}} \pi_i^{st.}(\mathbf{p},\mathbf{x}_i)\} = \pi_i^{st.+}(\mathbf{p}).$$

We note that this expression for $N_i(\mathbf{p},\mathbf{x}_i)$ can be obtained from (24) with

$$\gamma_i(\mathbf{p},\mathbf{x}_i) = (\delta_{\mathbf{u}_i,\mathbf{u}_i^*} - \delta_{\mathbf{x}_i,\mathbf{x}_i^*})[\pi_i^{st.+}(\mathbf{p}) - \pi_i^{st.*}(\mathbf{p})].$$

An alternative choice for a function $N_i(\mathbf{p},\mathbf{x}_i)$ originates from a set of $2^T$ redundant constraints parameterized by the vector of statuses $\mathbf{w}_i = (w_i^1,...,w_i^T)$ with $w_i^t \in \{0,1\}$, $t \in \{1,...,T\}$. Let $\rho_{\mathbf{w}_i}(\mathbf{u}_i) = -\prod_{t \in \{1,...T\}}(u_i^t)^{w_i^t}(1-u_i^t)^{(1-w_i^t)}$. Clearly, $\rho_i^{\mathbf{w}_i}(\mathbf{u}_i) \leq 0$, $\forall \mathbf{x}_i \in X_i$, $\forall \mathbf{w}_i \in \{0,1\}^T$. These functions have the following properties, $\forall \mathbf{u}_i, \mathbf{w}_i \in \{0,1\}^T$:

$$\rho_i^{\mathbf{w}_i}(\mathbf{u}_i) = -1 \leftrightarrow \mathbf{w}_i = \mathbf{u}_i; \quad \rho_i^{\mathbf{w}_i}(\mathbf{u}_i) = 0 \leftrightarrow \mathbf{w}_i \neq \mathbf{u}_i, \quad (30)$$

which implies that $\rho_i^{\mathbf{w}_i}(\mathbf{u}_i)$ can be expressed as $\rho_i^{\mathbf{w}_i}(\mathbf{u}_i) = -\delta_{\mathbf{u}_i,\mathbf{w}_i}$. Define $X_{\rho_i^{\mathbf{w}_i} \neq 0} = \{\mathbf{x}_i \mid \mathbf{x}_i \in X_i, \mathbf{u}_i = \mathbf{w}_i\}$. We note that if for some $\mathbf{w}_i \in \{0,1\}^T$ all the points with $\mathbf{u}_i = \mathbf{w}_i$ are infeasible in the generator private feasible set $X_i$ (for example, due to initial conditions combined with minimum up/down time constraints), then $X_{\rho_i^{\mathbf{w}_i} \neq 0} = \emptyset$. Let us denote as $\overline{W}_i$ the set of all $\mathbf{w}_i \in \{0,1\}^T$ with nonempty $X_{\rho_i^{\mathbf{w}_i} \neq 0}$. For $\forall \mathbf{w}_i \in \overline{W}_i$, Proposition 7 yields

$$\mu_i^{+\mathbf{w}_i \max} = \min_{\mathbf{x}_i \in X_{\rho_i^{\mathbf{w}_i} \neq 0}}[\pi_i^{st.+}(\mathbf{p}) - \pi_i^{st.}(\mathbf{p},\mathbf{x}_i)] = \pi_i^{st.+}(\mathbf{p}) - \max_{\substack{\mathbf{x}_i, \\ s.t. \\ \mathbf{x}_i \in X, \\ \mathbf{u}_i = \mathbf{w}_i}} \pi_i^{st.}(\mathbf{p},\mathbf{x}_i).$$

Using (28), we arrive at

$$\mu_i^{+\mathbf{w}_i \max} = \pi_i^{st.+}(\mathbf{p}) - \pi_i^{st.\max}(\mathbf{p},\mathbf{w}_i), \quad \forall \mathbf{w}_i \in \overline{W}_i. \quad (31)$$

Obviously, we have $\mu_i^{+\mathbf{w}_i \max} \geq 0$, $\forall \mathbf{w}_i \in W_i$. From (30) it follows that $\rho_i^{\mathbf{w}_i}(\mathbf{u}_i)\rho_i^{\mathbf{w}_i'}(\mathbf{u}_i) = 0$, $\forall \mathbf{u}_i, \mathbf{w}_i, \mathbf{w}_i' \in \{0,1\}^T$, $\mathbf{w}_i \neq \mathbf{w}_i'$, which entails $X_{\rho_i^{\mathbf{w}_i} \neq 0} \cap X_{\rho_i^{\mathbf{w}_i'} \neq 0} = \emptyset$, $\forall \mathbf{w}_i, \mathbf{w}_i' \in \{0,1\}^T$, $\mathbf{w}_i \neq \mathbf{w}_i'$. Therefore, the conditions of Proposition 9 are met and we have $M_i^+(\mathbf{p}) = (\times_{\mathbf{w}_i \in \overline{W}_i}[0, \mu_i^{+\mathbf{w}_i \max}]) \times R_{\geq 0}^{2T - |\overline{W}_i|}$. Hence, $\forall \mathbf{w}_i \in \overline{W}_i$, the functions $\rho_i^{\mathbf{w}_i}(\mathbf{u}_i)$ with the associated non-negative multipliers $\mu_i^{+\mathbf{w}_i \max}$ satisfy (16) and (18). It is straightforward to check that (17) is satisfied as well. Indeed, we have

$$-\sum_{\mathbf{w}_i \in \overline{W}_i} \mu_i^{+\mathbf{w}_i \max} \rho_i^{\mathbf{w}_i}(\mathbf{u}_i^*) = \sum_{\mathbf{w}_i \in \overline{W}_i} \mu_i^{+\mathbf{w}_i \max} \delta_{\mathbf{u}_i^*, \mathbf{w}_i} = \mu_i^{+\mathbf{u}_i^* \max} = \pi_i^{st.+}(\mathbf{p}) - \pi_i^{st.\max}(\mathbf{p}, \mathbf{u}_i^*),$$

which, using (29), is transformed to $-\sum_{\mathbf{w}_i \in \overline{W}} \mu_i^{+\mathbf{w}_i \max} \rho_i^{\mathbf{w}_i}(\mathbf{u}_i^*) = \pi_i^{st.+}(\mathbf{p}) - \pi_i^{st.*}(\mathbf{p})$. Therefore, (17) holds, and Proposition 4 entails that the introduction of the redundant constraints $\rho_i^{\mathbf{w}_i}(\mathbf{u}_i) \leq 0$, $\mathbf{w}_i \in \overline{W}_i$, with $\rho_i^{\mathbf{w}_i}(\mathbf{u}_i) = -\delta_{\mathbf{u}_i, \mathbf{w}_i}$, and the associated multipliers $\mu_i^{+\mathbf{w}_i \max}$ given by (31) results in zero uplift payment. (We note that only the constraint $\rho_i^{\mathbf{w}_i}(\mathbf{u}_i)$ with $\mathbf{w}_i = \mathbf{u}_i^*$ contributes to the uplift reduction, while the rest of $\rho_i^{\mathbf{w}_i}(\mathbf{u}_i)$ make no contribution to the uplift payment.) These constraints result in

$$N_i(\mathbf{p}, \mathbf{x}_i) = -\sum_{\mathbf{w}_i \in \overline{W}} \mu_i^{+\mathbf{w}_i \max} \rho_i^{\mathbf{w}_i}(\mathbf{u}_i) = \pi_i^{st.+}(\mathbf{p}) - \pi_i^{st.\max}(\mathbf{p}, \mathbf{u}_i), \quad (32)$$

and the amended profit function has the form $\pi_i(\mathbf{p}, \mathbf{x}_i) = \pi_i^{st.}(\mathbf{p}, \mathbf{x}_i) + \pi_i^{st.+}(\mathbf{p}) - \pi_i^{st.\max}(\mathbf{p}, \mathbf{u}_i)$. It is straightforward to check that the expression (32) for $N_i(\mathbf{p}, \mathbf{x}_i)$ can be expressed in terms of (24) with

$$\gamma_i(\mathbf{p}, \mathbf{x}_i) = \pi_i^{st.+}(\mathbf{p}) - \pi_i^{st.\max}(\mathbf{p}, \mathbf{u}_i) - \delta_{\mathbf{x}_i, \mathbf{x}_i^*}[\pi_i^{st.+}(\mathbf{p}) - \pi_i^{st.*}(\mathbf{p})].$$

We note that for the given market price $\mathbf{p}$ and statuses of the unit, the new terms in the profit function are constant. In the case of one-period market model with no intertemporal constraints (ramp, minimum up/down time constraints, etc.), we have $\delta_{u_i, 0} = (1 - u_i)$, $\delta_{u_i, 1} = u_i$, which entails

$$\pi_i^{st.\max}(p, u_i) = (1 - u_i)\pi_i^{st.\max}(p, 0) + u_i \pi_i^{st.\max}(p, 1) = u_i \pi_i^{st.\max}(p, 1),$$

where we used $\pi_i^{st.\max}(p, 0) = 0$. This allows expressing (32) as $N_i(p, x_i) = \pi_i^{st.+}(p) - u_i \pi_i^{st.\max}(p, 1)$. Thus, if the generator has the lost profit (i.e. it is offline in the centralized dispatch solution, but operation at the given market price $p$ would result in the profit $\pi_i^{st.+}(p) > 0$), then $N_i(p, x_i)$ compensates the generator for the lost profit $\pi_i^{st.+}(p)$ if it complies with the centralized dispatch solution. Likewise, if the generator operates at a loss (i.e. the generator is online in the centralized dispatch solution, but it is not recovering its cost at the given market price $p$ and, hence, would prefer to be offline), then $\pi_i^{st.+}(p) = 0$, $\pi_i^{st.\max}(p, 1) < 0$. In this case, $N_i(p, x_i)$ ensures that the generator receives zero profit and, as a result, fully recovers it cost if it follows the centralized dispatch.

## VII. Application to a producer with constant marginal cost of output in a single-period power market

Consider a producer operating a generating unit without intertemporal constraints in a single-period power market. The producer is assumed to have the constant marginal cost of output $a$, start-up cost $w$, minimum/maximum capacity limit $g^{\min} / g^{\max}$ with $0 < g^{\min} < g^{\max}$. Thus, the producer cost function has the form $C(x) = wu + ag$ defined on the producer private feasible set $X = \{(u,g) \mid u \in \{0,1\}, g \in R, ug^{\min} \leq g \leq ug^{\max}\}$. We consider the standard revenue function of the form $R^{st.}(p,x) = pg$ with the market price $p$, which is considered as being fixed by some pricing principle. Initially, the unit is assumed to be offline. Let $x^* = (u^*, g^*)$ denote the value of the producer status variable/output volume according to the centralized dispatch solution. We show that to construct the amendment function $N(p,x)$ that satisfies (2) - (4) it is sufficient to consider the following redundant constraints: $ug^{\min} - g \leq 0$, $g - ug^{\max} \leq 0$, $u - 1 \leq 0$. (We note that this list of the redundant constraints is not exhaustive since there are the redundant constraints that are not expressed as some linear combination of these constraints, e.g. $\max[ug^{\min} - g, g - ug^{\max}] \leq 0$.) Thus, we have

$$N(p,x) = \mu^1(g - ug^{\min}) + \mu^2(ug^{\max} - g) + \mu^3(1-u). \quad (33)$$

If $g^* = 0$, then for $p < a + w/g^{\max}$ the uplift payment is zero and we set $\mu^1 = \mu^2 = \mu^3 = 0$. For $p \geq a + w/g^{\max}$, from (13) we have $\mu^1 = 0$. Also, the constraint $g - ug^{\max} \leq 0$ is satisfied as equality at $g^* = 0$ and, therefore, does not contribute to the uplift. Consequently, we set $\mu^2 = 0$. Thus, we are left with $u - 1 \leq 0$. Proposition 7 gives $\mu^{3\max} = \pi^{st.+}(p)$. It is straightforward to check that the redundant constraint and the multiplier satisfy (16) – (18), which results in $N(p,x) = \pi^{st.+}(p)(1-u)$. In this case, the amended profit function is $\pi(p,x) = (p-a)g - wu + \pi^{st.+}(p)(1-u)$.

If $g^* = g^{\min}$, then we set $\mu^1 = \mu^3 = 0$ as the corresponding constraints do not affect the uplift at $x^*$. Proposition 7 entails $\mu^{2\max} = [\pi^{st.+}(p) - \pi^{st.*}(p)]/(g^{\max} - g^*)$, which gives $N(p,x) = [\pi^{st.+}(p) - \pi^{st.*}(p)](ug^{\max} - g)/(g^{\max} - g^*)$ with the redundant constraint and the multiplier satisfying (16) – (18). The amended profit function is given by

$$\pi(p,x) = (p-a)g - wu - [\pi^{st.+}(p) - \pi^{st.*}(p)]g/(g^{\max} - g^*) + [\pi^{st.+}(p) - \pi^{st.*}(p)]g^{\max}u/(g^{\max} - g^*).$$

If $g^{\min} < g^* < g^{\max}$, then for $p < a + w/g^{\max}$ we have $\mu^3 = 0$ as a consequence of (13). The conditions (17) – (18) give $\mu^1 = -\pi^{st.}(p,1,g^{\max})/(g^{\max} - g^{\min})$, $\mu^2 = -\pi^{st.}(p,1,g^{\min})/(g^{\max} - g^{\min})$. The resulting expression for the amendment function is $N(p,x) = -\pi^{st.}(p,x)$, which corresponds to the redundant constraint $(p-a)g - wu \leq 0$ and yields the identically zero amended profit function on $X$. In the case of $p \geq a + w/g^{\max}$, from (13) we have $\mu^1 = 0$, while the constraint $u - 1 \leq 0$ is satisfied as equality at $x^*$ and makes no contribution to the uplift. Therefore, we set $\mu^3 = 0$. Proposition 7 gives $\mu^{2\max} = (p-a)$. These redundant constraints and the corresponding multipliers satisfy (16) – (18) and give $N(p,x) = (p-a)(ug^{\max} - g)$, which produces the amended profit function $\pi(p,x) = \pi^{st.+}(p)u$.

Likewise, if $g^* = g^{\max}$, then for $p \geq a + w/g^{\max}$ the uplift is zero and no revenue function amendment is needed. For $p < a + w/g^{\max}$, we have $\mu^3 = 0$ from (13) and we also set $\mu^2 = 0$ since $g - ug^{\max} \leq 0$ is satisfied as equality at $x^*$. Proposition 7 gives $\mu^{1\max} = -\pi^{st.*}(p)/(g^{\max} - g^{\min})$ with $\pi^{st.*}(p) = (p-a)g^{\max} - w$ resulting in $N(p,x) = -\pi^{st.*}(p)(g - ug^{min})/(g^{max} - g^{min})$. It is straightforward to verify that the conditions (16) – (18) hold, and the amended profit function is expressed as $\pi(p,x) = (p-a)g - wu - \pi^{st.*}(p)(g - g^{min}u)/(g^{max} - g^{min})$.

Let us construct the profit function amendments resulting from the application of (27). We define the extended value cost function $\bar{C}(x)$ on $Conv\{X\}$ as $\bar{C}(x) = C(x)$ if $x \in X$, and $\bar{C}(x) = +\infty$ if $x \in Conv\{X\} \setminus X$. We have $Conv\{X\} = \{(u,g) \mid u \in [0,1], g \in [0, g^{max}], g/g^{min} \leq u \leq g/g^{max}\}$. Let us introduce a function $f(x) = conv_{Conv\{X\}}\{\bar{C}(x) - \delta_{x,x^*}[\pi^{st.+}(p) - \pi^{st.*}(p)]\}$ on $x \in Conv\{X\}$. For $x \in Conv\{X\}$, the function $f(x)$ defines a surface in $R^3$ with coordinates $(u, g, f)$. It can be shown that (27) yields the same expressions for $N(p,x)$ as above except for the case $g^{min} < g^* < g^{max}$. In this instance, for $x \in Conv\{X\}$ the function $f(x)$ corresponds to the highest of two planes in $R^3$ defined below. We have $f(x) = \max[f_1(x), f_2(x)]$ with $f_1(x)$ corresponding to the plane that contains the points $(0,0,0)$, $(1, g^{min}, w + ag^{min})$, $(1, g^*, w + ag^* - [\pi^{st.+}(p) - \pi^{st.*}(p)])$ and $f_2(x)$ corresponding to the plane that contains the points $(0,0,0)$, $(1, g^*, w + ag^* - [\pi^{st.+}(p) - \pi^{st.*}(p)])$, $(1, g^{max}, w + ag^{max})$. The straightforward computation gives

$$f_1(x) = wu + ag + [\pi^{st.+}(p) - \pi^{st.*}(p)](ug^{min} - g)/(g^* - g^{min}),$$
$$f_2(x) = wu + ag + [\pi^{st.+}(p) - \pi^{st.*}(p)](g - ug^{max})/(g^{max} - g^*).$$

This entails
$$N(p,x) = [\pi^{st.+}(p) - \pi^{st.*}(p)]\min[(g - ug^{min})/(g^* - g^{min}); (ug^{max} - g)/(g^{max} - g^*)], \quad (34)$$
which can be obtained from the redundant constraint $\max[(ug^{min} - g)/(g^* - g^{min}); (g - ug^{max})/(g^{max} - g^*)] \leq 0$ with the multiplier $\mu = \pi^{st.+}(p) - \pi^{st.*}(p)$. For $p \geq a + w/g^{max}$, the resulting amended profit function has the form $\pi(p,x) = u\pi^{st.+}(p) + \min[g - ug^*; 0](p-a)(g^{max} - g^{min})/(g^* - g^{min})$. For $p < a + w/g^{max}$, we have $\pi(p,x) = \min[(ug^* - g)\pi^{st.}(1, p, g^{min})/(g^* - g^{min}); (g - ug^*)\pi^{st.}(1, p, g^{max})/(g^{max} - g^*)]$.

It is illustrative to repeat the analysis when the status variable $u$ is expressed in terms of the output volume $g$ as $u = \theta(g)$ and is excluded from the consideration. In this case, the private feasible set of the generator is expressed as $G = \{g \mid g \in \{0\} \cup [g^{min}, g^{max}]\}$, while the cost function is $C(p,g) = ag + w\theta(g)$. We also have $Conv\{G\} = \{g \mid 0 \leq g \leq g^{max}\}$. The analysis above implies that to construct the amendment function $N(p,g)$ that yields zero uplift payment it is sufficient to consider the redundant constraints $\theta(g)g^{min} - g \leq 0$, $g - \theta(g)g^{max} \leq 0$, $\theta(g) - 1 \leq 0$. Such an approach gives $N(p,x) = \mu^1[g - \theta(g)g^{min}] + \mu^2[\theta(g)g^{max} - g] + \mu^3[1 - \theta(g)]$ with the same values of the multipliers $\mu^1$, $\mu^2$, $\mu^3$ as in (33).

Application of (27) results in the following expressions for the functions $N(p,g)$ and $\pi(p,g)$.

If $g^* = 0$, then for $p < a + w/g^{max}$ the uplift payment is not needed. For $p \geq a + w/g^{max}$, we have $N(p,g) = \pi^{st.+}(p) - \pi^{st.}(p,g)$, which can be obtained from the redundant constraint $\pi^{st.}(p,g) \leq \pi^{st.+}(p)$. The resulting amended profit function has the form $\pi(p,g) = \pi^{st.+}(p)$ and is constant.

If $g^* = g^{min}$, then for $p \geq a + w/g^{max}$ we have $N(p,g) = \pi^{st.+}(p)\theta(g) - \pi^{st.}(p,g)$, which is associated with the redundant constraint $\pi^{st.}(p,g) \leq \pi^{st.+}(p)\theta(g)$. For the amended profit function we have $\pi(p,g) = \pi^{st.+}(p)\theta(g)$. In the case of $p < a + w/g^{max}$, the calculation gives $N(p,g) = -\theta(g)(g^{max} - g)\pi^{st.}(p,g^{min})/(g^{max} - g^{min})$, which can be obtained from the redundant constraint $g \leq \theta(g)g^{max}$. (We note that the presence of $\theta(g)$ is critical in this constraint since the redundant constraint $g \leq g^{max}$ makes no contribution to the uplift payment as the associated multiplier is zero due to (13).) In this case, we have $\pi(p,g) = \min\{(g - g^{min})\pi^{st.}(p,g^{max})/(g^{max} - g^{min});0\}$.

If $g^{min} < g^* < g^{max}$, then in the case of $p \geq a + w/g^{max}$ we have $N(p,g) = \min[\{g[\pi^{st.+}(p) - \pi^{st.*}(p)] + w[\theta(g)g^* - g]\}/g^*; \pi^{st.+}(p) - \pi^{st.}(p,g)]$, which corresponds to the redundant constraint $\max[\{g[\pi^{st.*}(p) - \pi^{st.+}(p)] + w[g - \theta(g)g^*]\}/g^*; \pi^{st.}(p,g) - \pi^{st.+}(p)] \leq 0$. This constraint is equivalent to the set of the redundant constraints $\{g[\pi^{st.*}(p) - \pi^{st.+}(p)] + w[g - \theta(g)g^*]\}/g^* \leq 0$, which can be transformed to $g - g^{max} \leq 0$, and $\pi^{st.}(p,g) - \pi^{st.+}(p) \leq 0$. (However, as we noted in Section IV, the equivalent set of constraints may result in the different amendment function.) The amended profit function is $\pi(p,g) = \pi^{st.+}(p)\min[g;g^*]/g^*$. For $p < a + w/g^{max}$, we have

$$N(p,g) = \min[-\pi^{st.}(p,g); \{g^*[\pi^{st.}(p,g) - \pi^{st.}(p,g^{max})] + w[\theta(g)g^{max} - g]\}/(g^{max} - g^*)], \quad (35)$$

which is associated with the redundant constraint $-N(p,g) \leq 0$. (We note that the constraint $-N(p,g) \leq 0$ is equivalent to the set of the redundant constraints $\pi^{st.}(p,g) \leq 0$ and $g^*[\pi^{st.}(p,g^{max}) - \pi^{st.}(p,g)] + w[g - \theta(g)g^{max}] \leq 0$.) The amended profit function is given by $\pi(p,g) = \min[\pi^{st.}(p,g^{max})(g - g^*)/(g^{max} - g^*);0]$.

If $g^* = g^{max}$, then for $p \geq a + w/g^{max}$ the uplift payment is not needed. For $p < a + w/g^{max}$, we have $N(p,g) = -\pi^{st.}(p,g)$, which corresponds to the redundant constraint $\pi^{st.}(p,g) \leq 0$. The resulting amended profit function is identically zero.

## VIII. Numerical example

In this section we apply the amended profit function expressions (27) and (33) for the Scarf example [5] (adapted according to [22]), which describes the uninode single-period power market with fixed demand and three types of the power plants ("Smokestack", "High Tech", and "Med Tech") with the constant marginal costs of output. The unit parameters are given below.

**Table 1. Characteristics of the generating units**

|  | Minimum capacity limit, MW | Maximum capacity limit, MW | Marginal cost of output, $/MW | Start-up cost, $ |
|---|---|---|---|---|
| Smokestack | 0 | 16 | 3 | 53 |
| High Tech | 0 | 7 | 2 | 30 |
| Med Tech | 2 | 6 | 7 | 0 |

It is assumed that the power system has 6, 5, and 5 units of each type, respectively. Initially, all the units are offline. We consider two scenarios with demand equal 10 MWh and 40 MWh, respectively, with the standard revenue function of the form $R^{st.}(p,x) = pg$. We apply the convex hull pricing mechanism [20]-[22] to set the market price. The centralized dispatch outcomes and the market prices are given below, [22].

**Table 2. Dispatch and pricing outcome**

| Demand, MWh | Smokestack | | High Tech | | Med Tech | | Market price $p_{CHP}$, $/MWh | Total uplift, $ |
|---|---|---|---|---|---|---|---|---|
| | Number of units online | Each unit output, MWh | Number of units online | Each unit output, MWh | Number of units online | Each unit output, MWh | | |
| 10 | 0 | 0 | 1 | 7 | 1 | 3 | 6.2857 | 2.143 |
| 40 | 1 | 16 | 3 | 7 | 1 | 3 | 6.3125 | 2.438 |

If demand is 10 MWh, then the online High Tech unit sets the market price, and the uplift payment is made only to the operating Med Tech unit to compensate its output cost. Thus, $u_{MT}^* = 1$, $g_{MT}^* = 3$, $\pi_{MT}^{st.*}(p_{CHP}) = -\$2.143$, $u_{MT}^{st.+} = g_{MT}^{st.+} = 0$, $\pi_{MT}^{st.+}(p_{CHP}) = \$0$.

The application of (27) gives the following. From (34) we obtain the amendment function
$$N_{MT}(p_{CHP}, x_{MT}) = 2.143 \min[g_{MT} - 2u_{MT}; 2u_{MT} - g_{MT}/3],$$
which reaches its maximum value (equal the uplift payment of \$2.143) at $u_{MT}^* = 1$, $g_{MT}^* = 3$, and vanishes at $u_{MT} = g_{MT} = 0$, which is the optimal point of $\pi_{MT}^{st.}(p_{CHP}, x_{MT})$. It is straightforward to check that $N_{MT}(p_{CHP}, x_{MT}) \geq 0$, $\forall x_{MT} \in X_{MT}$, and such a choice for $N_{MT}(p_{CHP}, x_{MT})$ can be realized by the introduction of the redundant constraint $\min[g_{MT} - 2u_{MT}; 2u_{MT} - g_{MT}/3] \geq 0$ with the associated multiplier $\mu = 2.143$. The resulting amended profit function is
$$\pi_{MT}(p_{CHP}, x_{MT}) = 1.429 \min[g_{MT} - 3u_{MT}; 3u_{MT} - g_{MT}] = -1.429 |g_{MT} - 3u_{MT}|.$$

It is illustrative to consider the outcome of (27) if the private feasible set of a Med Tech unit is described using the output variable only. From (35) we obtain $N_{MT}(p_{CHP}, g_{MT}) = 0.714 \min[g_{MT}; 6 - g_{MT}]$. This expression for $N_{MT}(p_{CHP}, g_{MT})$ corresponds to the redundant constraint $\min[g_{MT}; 6 - g_{MT}] \geq 0$, which is equivalent to the set of the redundant constraints $0 \leq g_{MT}$, $g_{MT} \leq 6$. The resulting amended profit function is given by $\pi_{MT}(p_{CHP}, g_{MT}) = 2.143 \min[3 - g_{MT}; 0]$.

The alternative amendment function can be obtained from (33), which yields the redundant constraints $2u_{MT} - g_{MT} \leq 0$ and $g_{MT} - 6u_{MT} \leq 0$ with the multipliers equal 1.071 and 0.357, respectively. The resulting amendment function equals $-\pi_{MT}^{st.}(p_{CHP}, x_{MT})$ and produces identically zero amended profit function.

Now, we consider a scenario with demand equal 40 MWh. Without the uplift payment, the online Med Tech unit has a loss of \$2.063 and each of the two offline High Tech units has the lost profit in the amount of \$0.188. First, let us consider the online Med Tech. The application of (27) for the formulation involving the private feasible set $X_{MT}$ is given by (34), which yields $N_{MT}(p_{CHP}, x_{MT}) = 2.063 \min[g_{MT} - 2u_{MT}; 2u_{MT} - g_{MT}/3]$. (This expression for the amendment function $N_{MT}(p_{CHP}, x_{MT})$ can be obtained from the redundant constraint $\min[g_{MT} - 2u_{MT}; 2u_{MT} - g_{MT}/3] \geq 0$.) The corresponding amended profit function is expressed as $\pi_{MT}(p_{CHP}, x_{MT}) = -1.375 |g_{MT} - 3u_{MT}|$. In case of the private feasible set $G_{MT}$, (27) gives $N_{MT}(p_{CHP}, g_{MT}) = 0.688 \min[g_{MT}; 6 - g_{MT}]$. We note that this expression for $N_{MT}(p_{CHP}, g_{MT})$ can be obtained from the redundant constraint $\min[g_{MT}; 6 - g_{MT}] \geq 0$, which is equivalent to the redundant constraints $0 \leq g_{MT}$, $g_{MT} \leq 6$ belonging to the private feasible set of the Med Tech unit. The amended profit function is expressed as $\pi_{MT}(p_{CHP}, g_{MT}) = \min[0; 4.126 - 1.375 g_{MT}]$. Another choice for the amendment function is given by (33), which entails $N_{MT}(p_{CHP}, x_{MT}) = 1.031(g_{MT} - 2u_{MT}) + 0.344(6u_{MT} - g_{MT})$ with identically zero amended profit function.

For the High Tech unit, the application of both (27) and (33) for the producer private feasible set formulated as $X_{HT}$ yields $N_{HT}(p_{CHP}, x_{HT}) = 0,188(1 - u_{HT})$ and $\pi_{HT}(p_{CHP}, x_{HT}) = 0,188 + 4.313 g_{HT} - 30.188 u_{HT}$. An alternative choice for the amendment function is given by (27) for the case of the private feasible set $G_{HT}$. We have $N_{HT}(p_{CHP}, g_{HT}) = 0.188 + 30\theta(g_{HT}) - 4.313 g_{HT}$, which can be deduced from the redundant constraint $\pi_{HT}^{st.}(p_{CHP}, g_{HT}) - \pi_{HT}^{st.+}(p_{CHP}) \leq 0$. The resulting amended profit function is constant on $G_{HT}$ and equals 0.188.

## IX. Conclusion

We considered the possibly non-linear redundant constraints that stay redundant if the power balance constraint is excluded from the constraint set of the centralized dispatch optimization problem and showed that the introduction of this type of the redundant constraints leaves the duality gap unaffected. For each producer, we studied the redundant constraints that are satisfied on its private feasible set and proved that any set of the redundant constraints, which belong to this special type, corresponds to the producer revenue function amendment that is non-negative (on the producer private feasible set) and leaves the maximum profit of the producer unaffected. Likewise, for any such amendment function, one can indicate (generally non-unique) set of the redundant constraints. Consequently, the uplift payment is potentially lowered by the introduction of these constraints.

We studied the properties of the redundant constraints and formulated necessary and sufficient conditions for a given set of the redundant constraints and the associated multipliers to yield zero uplift payment. For each producer, we explicitly construct the general expression for the producer revenue function amendment that is non-negative on the producer private feasible set, leaves the maximum profit of the producer unaffected, and results in zero uplift payment for the market player. This allows identifying the family of the redundant constraints that corresponds to the revenue function amendment of the producer and yields zero uplift payment for this market player. In case of the uniform price for power, we constructed one universal redundant constraint (given by the sum of these properly rescaled individual redundant constraints) that could be introduced directly in the centralized dispatch optimization problem to yield zero total uplift payment after the Lagrangian relaxation procedure is applied to both this constraint and the power balance constraint.

Thus, in the case of a uniform market price, it suffices to introduce just one redundant constraint in the centralized dispatch optimization problem to eliminate all the uplift payments. If the uniform market price is set using the convex hull pricing method, the set of the market prices and each producer maximum profit are unaffected by the redundant constraint.

**Appendix**

In this section we study the properties of the redundant constraints and establish a necessary condition for a given set of the redundant constraints to yield zero uplift payment for the producer $i$. For $\forall l_i \in L_i$, consider an optimization problem

$$\pi_i^{st.+}(\mathbf{p}) = \min_{\mu_i^{l_i} \geq 0} \max_{\mathbf{x}_i \in X_i} \; [\pi_i^{st.}(\mathbf{p},\mathbf{x}_i) - \mu_i^{l_i} \rho_i^{l_i}(\mathbf{p},\mathbf{x}_i)]. \quad (36)$$

Let us define the corresponding set of minimizers $M_i^{+l_i}(\mathbf{p}) = \arg\min_{\mu_i^{l_i} \geq 0} \max_{\mathbf{x}_i \in X_i} \; [\pi_i^{st.}(\mathbf{p},\mathbf{x}_i) - \mu_i^{l_i} \rho_i^{l_i}(\mathbf{p},\mathbf{x}_i)]$. Clearly, $M_i^{+l_i}(\mathbf{p})$ is a closed convex set and $\{0\} \in M_i^{+l_i}(\mathbf{p})$. Since $M_i^{+l_i}(\mathbf{p}) \in R_{\geq 0}$, we have the following three possibilities for the set $M_i^{+l_i}(\mathbf{p})$: $M_i^{+l_i}(\mathbf{p}) = \{0\}$, or $M_i^{+l_i}(\mathbf{p}) = [0,a]$ with some $a \in R_{>0}$, or $M_i^{+l_i}(\mathbf{p}) = R_{\geq 0}$. Obviously, if $\rho_i^{l_i}(\mathbf{p},\mathbf{x}_i) < 0$, $\forall \mathbf{x}_i \in X_i$, then $M_i^{+l_i}(\mathbf{p}) = \{0\}$, while $\rho_i^{l_i}(\mathbf{p},\mathbf{x}_i) = 0$, $\forall \mathbf{x}_i \in X_i$, entails $M_i^{+l_i}(\mathbf{p}) = R_{\geq 0}$. Let us denote by $\bar{L}_i$ the subset of $L_i$ with bounded $M_i^{+l_i}(\mathbf{p})$. Thus, $\forall l_i \in \bar{L}_i$ we have $\rho_i^{l_i}(\mathbf{p},\mathbf{x}_i) \leq 0$, $\forall \mathbf{x}_i \in X_i$, and $\rho_i^{l_i}(\mathbf{p},\mathbf{x}'_i) \neq 0$ for some $\mathbf{x}'_i \in X_i$, which means that a set $X_{\rho_i^{l_i} \neq 0} = \{\mathbf{x}_i \mid \mathbf{x}_i \in X_i, \rho_i^{l_i}(\mathbf{p},\mathbf{x}_i) \neq 0\}$ is nonempty. For $\forall l_i \in \bar{L}_i$, let $\mu_i^{+l_i \max}$ denote the maximum element of $M_i^{+l_i}(\mathbf{p})$. The following statement gives a straightforward way to calculate $\mu_i^{+l_i \max}$ under some simplifying technical assumption.

*Proposition 7.* If for $l_i \in \bar{L}_i$ the function $[\pi_i^{st.}(\mathbf{p},\mathbf{x}_i) - \pi_i^{st.+}(\mathbf{p})]/\rho_i^{l_i}(\mathbf{p},\mathbf{x}_i)$ has a minimum value on $X_{\rho_i^{l_i} \neq 0}$, then

$$\mu_i^{+l_i \max} = \min_{\mathbf{x}_i \in X_{\rho_i^{l_i} \neq 0}} \; [\pi_i^{st.}(\mathbf{p},\mathbf{x}_i) - \pi_i^{st.+}(\mathbf{p})]/\rho_i^{l_i}(\mathbf{p},\mathbf{x}_i). \quad (37)$$

Proof. Let $\mu_i^{+l_i} \in M_i^{+l_i}(\mathbf{p})$, then (36) implies $\mu_i^{+l_i} \rho_i^{l_i}(\mathbf{p},\mathbf{x}_i) \geq \pi_i^{st.}(\mathbf{p},\mathbf{x}_i) - \pi_i^{st.+}(\mathbf{p})$, $\forall \mathbf{x}_i \in X_i$, which entails $\mu_i^{+l_i} \leq \min_{\mathbf{x}_i \in X_{\rho_i^{l_i} \neq 0}} \; [\pi_i^{st.}(\mathbf{p},\mathbf{x}_i) - \pi_i^{st.+}(\mathbf{p})]/\rho_i^{l_i}(\mathbf{p},\mathbf{x}_i)$. Let $\bar{\mu}_i^{l_i} = \min_{\mathbf{x}_i \in X_{\rho_i^{l_i} \neq 0}} \; [\pi_i^{st.}(\mathbf{p},\mathbf{x}_i) - \pi_i^{st.+}(\mathbf{p})]/\rho_i^{l_i}(\mathbf{p},\mathbf{x}_i)$. Thus, $\bar{\mu}_i^{l_i} \geq \mu_i^{+l_i}$, $\forall \mu_i^{+l_i} \in M_i^{+l_i}(\mathbf{p})$. Now, we show that $\bar{\mu}_i^{l_i} \in M_i^{+l_i}(\mathbf{p})$. For $\mathbf{x}'_i \in \arg\min_{\mathbf{x}_i \in X_{\rho_i^{l_i} \neq 0}} \; [\pi_i^{st.}(\mathbf{p},\mathbf{x}_i) - \pi_i^{st.+}(\mathbf{p})]/\rho_i^{l_i}(\mathbf{p},\mathbf{x}_i)$, we have $\pi_i^{st.+}(\mathbf{p}) = \pi_i^{st.}(\mathbf{p},\mathbf{x}'_i) - \bar{\mu}_i^{l_i} \rho_i^{l_i}(\mathbf{p},\mathbf{x}'_i)$. Since $\pi_i^{st.+}(\mathbf{p}) \geq \pi_i^{st.}(\mathbf{p},\mathbf{x}_i) - \bar{\mu}_i^{l_i} \rho_i^{l_i}(\mathbf{p},\mathbf{x}_i)$, $\forall \mathbf{x}_i \in X_i$, we conclude that $\pi_i^{st.+}(\mathbf{p}) = \max_{\mathbf{x}_i \in X_i} \; [\pi_i^{st.}(\mathbf{p},\mathbf{x}_i) - \mu_i^{l_i} \rho_i^{l_i}(\mathbf{p},\mathbf{x}_i)]$ and $\bar{\mu}_i^{l_i} \in M_i^{+l_i}(\mathbf{p})$. Therefore, $\mu_i^{+l_i \max} = \bar{\mu}_i^{l_i}$. Proposition is proved.

We note that if $\rho_i^{l_i}(\mathbf{p},\mathbf{x}_i^{st.+}) \neq 0$, $l_i \in \bar{L}_i$, then (37) yields $\mu_i^{+l_i \max} = 0$, which agrees with (13). Also, if $\rho_i^{l_i}(\mathbf{p},\mathbf{x}_i^*) \neq 0$, $l_i \in \bar{L}_i$, then (37) implies $\mu_i^{+l_i \max} \leq [\pi_i^{st.*}(\mathbf{p}) - \pi_i^{st.+}(\mathbf{p})]/\rho_i^{l_i}(\mathbf{p},\mathbf{x}_i^*)$, which entails $\mu_i^{+l_i \max} \rho_i^{l_i}(\mathbf{p},\mathbf{x}_i^*) \geq \pi_i^{st.*}(\mathbf{p}) - \pi_i^{st.+}(\mathbf{p})$. Therefore, if $\mu_i^{+l_i} \in M_i^{+l_i}(\mathbf{p})$, then $\mu_i^{+l_i}$ and $\rho_i^{l_i}(\mathbf{p},\mathbf{x}_i)$ satisfy (16) and (18). However, if (17) holds for $\mu_i^{+l_i} \in M_i^{+l_i}(\mathbf{p})$ and $\rho_i^{l_i}(\mathbf{p},\mathbf{x}_i)$, $l_i \in \bar{L}_i$, then $\mu_i^{+l_i} = \mu_i^{+l_i \max}$. The following statement gives a relation between $M_i^+(\mathbf{p})$ and the sets $M_i^{+l_i}(\mathbf{p})$.

*Proposition 8.* $M_i^+(\mathbf{p}) \subset \underset{l_i \in L_i}{\times} M_i^{+l_i}(\mathbf{p})$.

Proof: For a given $k \in L_i$, define a projection operation $P^k$ to a hyperplane in $R^{|L_i|}$ defined by $\mu_i^k = 0$: $P^k \mu = (\mu_i^1,...,\mu_i^{k-1}, 0, \mu_i^{k+1},...,\mu_i^{L_i})$. For $\forall \mu_i^+ \in M_i^+(\mathbf{p})$, we have

$$\pi_i^{st.}(\mathbf{p}, \mathbf{x}_i) \leq \pi_i^{st.}(\mathbf{p}, \mathbf{x}_i) - (P^k \mu_i^+)^T \rho_i(\mathbf{p}, \mathbf{x}_i) \leq \pi_i^{st.}(\mathbf{p}, \mathbf{x}_i) - \mu_i^{+T} \rho_i(\mathbf{p}, \mathbf{x}_i), \ \forall k \in L_i, \ \forall \mathbf{x}_i \in X_i.$$

Therefore, $\pi_i^{st.+}(\mathbf{p}) \leq \max_{\mathbf{x}_i \in X_i} [\pi_i^{st.}(\mathbf{p}, \mathbf{x}_i) - (P^k \mu_i^+)^T \rho_i(\mathbf{p}, \mathbf{x}_i)] \leq \pi_i^{st.+}(\mathbf{p})$. This implies $\pi_i^{st.+}(\mathbf{p}) = \max_{\mathbf{x}_i \in X_i} [\pi_i^{st.}(\mathbf{p}, \mathbf{x}_i) - (P^k \mu_i^+)^T \rho_i(\mathbf{p}, \mathbf{x}_i)]$. As a result, if $\mu_i^+ \in M_i^+(\mathbf{p})$, then $P^k \mu_i^+ \in M_i^+(\mathbf{p})$. The sequential application of the projection operation to various hyperplanes yields that if $\mu_i^+ = (\mu_i^{+1},...,\mu_i^{+|L_i|}) \in M_i^+(\mathbf{p})$, then $(0,...,0,\mu_i^{+l_i},0,..,0) \in M_i^+(\mathbf{p})$, which is equivalent to $\mu_i^{+l_i} \in M_i^{+l_i}(\mathbf{p})$. This implies $M_i^+(\mathbf{p}) \subset \underset{l_i \in L_i}{\times} M_i^{+l_i}(\mathbf{p})$. Proposition is proved.

Therefore, a set $M_i^+(\mathbf{p})$ is a subset of an $|L_i|$-dimensional box $\underset{l_i \in L_i}{\times} M_i^{+l_i}(\mathbf{p})$. We note that some facets of this box can be degenerate or unbounded if $M_i^{+l_i}(\mathbf{p}) = \{0\}$ or $M_i^{+l_i}(\mathbf{p}) \in R_{\geq 0}$ for some $l_i \in L_i$. Now we identify a special case when $M_i^+(\mathbf{p}) = \underset{l_i \in L_i}{\times} M_i^{+l_i}(\mathbf{p})$.

*Proposition 9.* If $X_{\rho_i^{l_i} \neq 0} \cap X_{\rho_i^{l_i'} \neq 0} = \emptyset$, $\forall l_i, l_i' \in L_i$, $l_i \neq l_i'$, then $M_i^+(\mathbf{p}) = \underset{l_i \in L_i}{\times} M_i^{+l_i}(\mathbf{p})$.

Proof: Proposition 8 entails that $M_i^+(\mathbf{p}) \subset \underset{l_i \in L_i}{\times} M_i^{+l_i}(\mathbf{p})$. Therefore, we need to show that $\underset{l_i \in L_i}{\times} M_i^{+l_i}(\mathbf{p}) \subset M_i^+(\mathbf{p})$. Let us choose an arbitrary $\mu_i^+ = (\mu_1^+,...,\mu_{|L_i|}^+) \in \underset{l_i \in L_i}{\times} M_i^{+l_i}(\mathbf{p})$ with $\mu_i^{+l_i} \in M_i^{+l_i}(\mathbf{p})$, $\forall l_i \in L_i$, and partition $X_i$ into the subsets with the elements belonging to one of $X_{\rho_i^{l_i} \neq 0}$, $l_i \in L_i$, and the rest of $X_i$. We have $X_i = [\underset{l_i \in L_i}{\cup} X_{\rho_i^{l_i} \neq 0}] \cup [X_i \setminus \underset{l_i \in L_i}{\cup} X_{\rho_i^{l_i} \neq 0}]$, therefore

$$\max_{\mathbf{x}_i \in X_i} [\pi_i^{st.}(\mathbf{p}, \mathbf{x}_i) - \mu_i^{+T} \rho_i(\mathbf{p}, \mathbf{x}_i)] =$$

$$\max \{\max_{l_i \in L_i} \max_{\mathbf{x}_i \in X_{\rho_i^{l_i} \neq 0}} [\pi_i^{st.}(\mathbf{p}, \mathbf{x}_i) - \mu_i^{+l_i} \rho_i^{l_i}(\mathbf{p}, \mathbf{x}_i)]; \max_{\mathbf{x}_i \in X_i \setminus \underset{l_i \in L_i}{\cup} X_{\rho_i^{l_i} \neq 0}} \pi_i^{st.}(\mathbf{p}, \mathbf{x}_i) \}.$$

Using $\max_{\mathbf{x}_i \in X_{\rho_i^{l_i} \neq 0}} [\pi_i^{st.}(\mathbf{p}, \mathbf{x}_i) - \mu_i^{+l_i} \rho_i^{l_i}(\mathbf{p}, \mathbf{x}_i)] \leq \max_{\mathbf{x}_i \in X_i} [\pi_i^{st.}(\mathbf{p}, \mathbf{x}_i) - \mu_i^{+l_i} \rho_i^{l_i}(\mathbf{p}, \mathbf{x}_i)] = \pi_i^{st.+}(\mathbf{p})$ and $\max_{\mathbf{x}_i \in X_i \setminus \underset{l_i \in L_i}{\cup} X_{\rho_i^{l_i} \neq 0}} \pi_i^{st.}(\mathbf{p}, \mathbf{x}_i) \leq \max_{\mathbf{x}_i \in X_i} \pi_i^{st.}(\mathbf{p}, \mathbf{x}_i) = \pi_i^{st.+}(\mathbf{p})$, we arrive at $\max_{\mathbf{x}_i \in X_i} [\pi_i^{st.}(\mathbf{p}, \mathbf{x}_i) - \mu_i^{+T} \rho_i(\mathbf{p}, \mathbf{x}_i)] \leq \pi_i^{st.+}(\mathbf{p})$. However, $\pi_i^{st.}(\mathbf{p}, \mathbf{x}_i) \leq \pi_i^{st.}(\mathbf{p}, \mathbf{x}_i) - \mu_i^{+T} \rho_i(\mathbf{p}, \mathbf{x}_i)$, $\forall \mathbf{x}_i \in X_i$. Therefore, $\pi_i^{st.+}(\mathbf{p}) \leq \max_{\mathbf{x}_i \in X_i} [\pi_i^{st.}(\mathbf{p}, \mathbf{x}_i) - \mu_i^{+T} \rho_i(\mathbf{p}, \mathbf{x}_i)]$. Consequently, $\pi_i^{st.+}(\mathbf{p}) = \max_{\mathbf{x}_i \in X_i} [\pi_i^{st.}(\mathbf{p}, \mathbf{x}_i) - \mu_i^{+T} \rho_i(\mathbf{p}, \mathbf{x}_i)]$, which yields $\mu_i^+ \in M_i^+(\mathbf{p})$. Since $\mu_i^+$ is an arbitrary element of $\underset{l_i \in L_i}{\times} M_i^{+l_i}(\mathbf{p})$, we conclude that $M_i^+(\mathbf{p}) = \underset{l_i \in L_i}{\times} M_i^{+l_i}(\mathbf{p})$. Proposition is proved.

We also note that if $X_{\rho_i^{l_i} \neq 0} \cap X_{\rho_i^{l_i'} \neq 0} = \emptyset$, $\forall l_i, l_i' \in L_i$, $l_i \neq l_i'$, then at most one of the redundant constraints contributes to the uplift (i.e. satisfies $\rho_i^{l_i}(\mathbf{p}, \mathbf{x}_i^*) \neq 0$) and the other redundant constraints from $L_i$ can be excluded from the consideration.

Motivated by an observation that if $\mu_i^{+l_i} \in M_i^{+l_i}(\mathbf{p})$, then $(0,...,0,\mu_i^{+l_i},0,...,0) \in M_i^+(\mathbf{p})$, we may indicate an nonempty subset of $M_i^+(\mathbf{p})$. For bounded $M_i^+(\mathbf{p})$, let us denote by $\Delta_i(\mathbf{p})$ the convex hull of the points $(0,...,0,\mu_i^{+l_i\max},0,...,0)$, $l_i \in L_i$, and $(0,...,0)$. In this case, since $M_i^+(\mathbf{p})$ is a convex set, we have $\Delta_i(\mathbf{p}) \subset M_i^+(\mathbf{p})$. For example, if all $\rho_i^{+l_i}(\mathbf{p},\mathbf{x}_i)$ are identical with nonempty $X_{\rho_i^{+l_i}\neq 0}$, then $M_i^+(\mathbf{p})$ is bounded and $\Delta_i(\mathbf{p}) = M_i^+(\mathbf{p})$.

Now we are ready to formulate the necessary condition that holds if a vector function $\rho_i(\mathbf{p},\mathbf{x}_i)$ produces zero uplift.

*Proposition 10.* Let $\rho_i(\mathbf{p},\mathbf{x}_i) \leq 0$, $\forall \mathbf{x}_i \in X_i$. If $\min_{\mu_i \in M_i^+(\mathbf{p})} U_i(\mathbf{p},\mu_i) = 0$, then $\sum_{l_i \in \bar{L}_i} \mu_i^{+l_i\max} \rho_i^{l_i}(\mathbf{p},\mathbf{x}_i^*) \leq \pi_i^{st.*}(\mathbf{p}) - \pi_i^{st.+}(\mathbf{p})$.

Proof. Due to $M_i^+(\mathbf{p}) \subset \times_{l_i \in L_i} M_i^{+l_i}(\mathbf{p})$, we have $\min_{\mu_i \in \times_{l_i \in L_i} M_i^{+l_i}(\mathbf{p})} \mu_i^T \rho_i(\mathbf{p},\mathbf{x}_i^*) \leq \min_{\mu_i \in M_i^+(\mathbf{p})} \mu_i^T \rho_i(\mathbf{p},\mathbf{x}_i^*)$.

Consequently, $\min_{\mu_i \in \times_{l_i \in L_i} M_i^{+l_i}(\mathbf{p})} \mu_i^T \rho_i(\mathbf{p},\mathbf{x}_i^*) = \sum_{l_i \in L_i} \min_{\mu_i^{l_i} \in M_i^{+l_i}(\mathbf{p})} \mu_i^{l_i} \rho_i^{l_i}(\mathbf{p},\mathbf{x}_i^*) = \sum_{l_i \in \bar{L}_i} \min_{\mu_i^{l_i} \in M_i^{+l_i}(\mathbf{p})} \mu_i^{l_i} \rho_i^{l_i}(\mathbf{p},\mathbf{x}_i^*)$, where we have used the fact that if $M_i^{+l_i}(\mathbf{p})$ is not bounded, then $\mu_i^{l_i} \rho_i^{l_i}(\mathbf{p},\mathbf{x}_i^*) = 0$. We also have $\mu_i^{+l_i\max} \rho_i^{l_i}(\mathbf{p},\mathbf{x}_i^*) \leq \min_{\mu_i^{l_i} \in M_i^{+l_i}(\mathbf{p})} \mu_i^{l_i} \rho_i^{l_i}(\mathbf{p},\mathbf{x}_i^*)$, $\forall l_i \in \bar{L}_i$. If $\min_{\mu_i \in M_i^+(\mathbf{p})} U_i(\mathbf{p},\mu_i) = 0$, then (15) implies $\min_{\mu_i \in M_i^+(\mathbf{p})} \mu_i^T \rho_i(\mathbf{p},\mathbf{x}_i^*) = \pi_i^{st.*}(\mathbf{p}) - \pi_i^{st.+}(\mathbf{p})$, which gives $\sum_{l_i \in \bar{L}_i} \mu_i^{+l_i\max} \rho_i^{+l_i}(\mathbf{p},\mathbf{x}_i^*) \leq \pi_i^{st.*}(\mathbf{p}) - \pi_i^{st.+}(\mathbf{p})$.

Proposition is proved.